\documentclass[12pt, leqno]{amsart}
\usepackage[english]{babel}
\usepackage{amssymb,bbm,xcolor,latexsym, amsmath}

\usepackage[leqno]{amsmath}
\catcode`\<=\active \def<{
\fontencoding{T1}\selectfont\symbol{60}\fontencoding{\encodingdefault}}
\newcommand{\assign}{:=}
\newcommand{\dueto}[1]{\textup{{(#1) }}}
\newcommand{\mathd}{\mathrm{d}}
\newcommand{\nocomma}{}
\newcommand{\noplus}{}

\newcommand{\tmem}[1]{{\em #1\/}}
\newcommand{\tmmathbf}[1]{\ensuremath{\boldsymbol{#1}}}
\newcommand{\tmmisc}[1]{\thanks{\textit{} #1}}
\newcommand{\tmams}[1]{\thanks{\textit{AMS 2010 classification:} #1}}
\newcommand{\tmname}[1]{\textsc{#1}}
\newcommand{\tmop}[1]{\ensuremath{\operatorname{#1}}}

\newcommand{\tmtextit}[1]{{\itshape{#1}}}
\newenvironment{proof*}[1]{\noindent\textbf{#1\ }}{\hspace*{\fill}$\Box$\medskip}
\newtheorem{lemma}{Lemma}
\newtheorem{proposition}{Proposition}
\newtheorem{theorem}{Theorem}
\newtheorem{definition}{Definition}
\newtheorem{question}{Question}
\begin{document}

\title{Differential subordination under change of law}
\tmmisc{Partially supported by ERC grant CHRiSHarMa 682402. The second author is a member of IUF }
\tmams{Primary  60G44; Secondary 60G46 }

\author{K.\;Domelevo}
\address{Insitut de Math{\'e}matiques de Toulouse.
Universit{\'e} Paul Sabatier, Toulouse, France}
\email{komla.domelevo@math.univ-toulouse.fr}

\author{S.\;Petermichl}
\address{Insitut de Math{\'e}matiques de Toulouse.
Universit{\'e} Paul Sabatier, Toulouse, France}
\email{stefanie.petermichl@math.univ-toulouse.fr}

{\maketitle}

\begin{abstract}
  We prove optimal $L^2$ bounds for a pair of Hilbert space valued
  differentially subordinate martingales under a change of law. The change of
  law is given by a process called a weight and sharpness in this context
  refers to the optimal growth with respect to the characteristic of the
  weight. The pair of martingales are adapted, uniformly integrable, and
  c{\`a}dl{\`a}g. Differential subordination is in the sense of Burkholder,
  defined through the use of the square bracket. In the scalar dyadic setting
  with underlying Lebesgue measure, this was proved by {\tmname{Wittwer}}
  {\cite{Wit2000a}}, where homogeneity was heavily used. Recent progress by
  {\tmname{Thiele}}--{\tmname{Treil}}--{\tmname{Volberg}}
  {\cite{ThiTreVol2015a}} and {\tmname{Lacey}} {\cite{Lac2015a}},
  independently, resolved the so--called non--homogenous case using discrete
  in time filtrations, where one martingale is a predictable multiplier of the
  other. The general case for continuous--in--time filtrations and pairs of
  martingales that are not necessarily predictable multipliers, remained open
  and is adressed here. As a very useful by--product, we give the explicit
  expression of a Bellman function of four variables for the weighted estimate
  of subordinate martingales with jumps. This construction includes an
  analysis of the regularity of this function as well as a very precise
  convexity needed to deal with the jump part.

  \ 
\end{abstract}

\section{Introduction}

\

The paper by {\tmname{Nazarov}}--{\tmname{Treil}}--{\tmname{Volberg}}
{\cite{NazTreVol1999a}} has set the groundwork for the early advances in
modern weighted theory in harmonic analysis and probability that started
around twenty years ago. In their paper the authors show necessary and
sufficient conditions for a dyadic martingale transform to be bounded in the
$L^2$ two--weight setting. The methodology of their proof could be used to get
the first sharp result in the real valued one-weight setting, for the dyadic
martingale transform {\cite{Wit2000a}}. Sharpness in this setting means best
control on growth with the necessary $A_2$ condition
\[ \sup_{\tau} \tmop{ess} . \sup_{\omega} \mathbbm{E} (w |
   \mathcal{F}_{\tau}) \mathbbm{E} (w^{- 1} | \mathcal{F}_{\tau}) \]
with $\tau$ adapted stopping times, where the dyadic filtration stands. Thus
this becomes
\[ \sup_I \left( \frac{1}{| I |} \int w \right) \left( \frac{1}{| I |} \int
   w^{- 1} \right) \]
where the supremum runs over all dyadic intervals.

\

The area of sharp weighted estimates has seen substantial progress with new,
beautiful proofs of {\tmname{Wittwer}}'s result and its extensions to the time
shifted martingales referred to as `dyadic shift'
{\cite{LacPetReg2010a}}{\cite{Tre2013a}}. Related, important \ questions in
harmonic analysis, such as boundedness of the Beurling--Ahlfors transform
{\cite{PetVol2002a}}, Hilbert transform {\cite{Pet2007a}}, general
Calderon-Zygmund operators {\cite{Hyt2012a}}{\cite{Ler2015a}}{\cite{Lac2015a}}
and beyond {\cite{BerFrePet2015a}}{\cite{Hyt2012a}} have been solved,
beautifully advancing profound understanding of the objects at hand.

During the early days of weighted theory in harmonic analysis, before optimal
weighted estimates were within reach, say, for the maximal operator or the
Hilbert transform {\cite{HunMucWhe1973a}} similar questions were asked in
probability theory, concerning stochastic processes with continuous in time
filtrations {\cite{BonLep1979a}}{\cite{IzuKaz1977a}}. The difficulty that
arises in the non--homogenous setting, typically seen when these processes
have jumps, were already observed back then and this restriction was made in
one form or another in these papers. Certain basic facts about weights do not
hold true for jump processes, such as the classical self improvement of the
$A_2$ characteristic of the weight {\cite{BonLep1979a}}. Another obstacle
typical for working with weights is the non-convexity of the set inspired by
the $A_2$ characteristic: $\{ r, s \in \mathbbm{R}_+ : 1 \leqslant r s
\leqslant Q \}$ with $Q > 1$. Such continuity--in--space assumptions still
appear regularly for these or other reasons when adressing weights, see
{\cite{BanOse2016a}}{\cite{Ose2016a}}.

{\tmname{Wittwer}}'s proof also uses the homogeneity that arises from the
dyadic filtration where the underlying measure is Lebesgue in a subtle but
crucial way. This homogeneity assumption has only recently been removed in the
papers {\cite{ThiTreVol2015a}} and {\cite{Lac2015a}}. These authors work with
discrete in time general filtrations with arbitrary underlying measure, where
one martingale is a predictable multiplier of the other. A direct passage
using the results for discrete in time filtrations to the continuous in time
filtration case where one uses Burkholder's definition
\begin{eqnarray}
   \label{defsub}&& Y \tmop{differentially} \tmop{subordinate} \tmop{to} X \\
 \nonumber &:\Leftrightarrow  & [X, X]_t - [Y, Y]_t \tmop{nonnegative} \tmop{and} \tmop{nondecreasing}
\end{eqnarray}
is only possible in very special cases, such as predictable multipliers of
stochastic integrals - this passage is explained in one of Burkholder's early
works on $L^p$ estimates for pairs of differentially subordinate martingales
{\cite{Bur1984a}}. (In full generality, this unweighted $L^p$ problem was only
much later resolved in {\cite{Wan1995a}}.)

In this article, we tackle the sharp weighted estimate in full generality,
using the notion of differentail subordination of Burkholder (\ref{defsub})
and the martingale $A_2$ characteristic,
\[ Q^{\mathcal{F}}_2 [w] = \sup_{\tau} \tmop{ess} . \sup_{\omega}  (w)_{\tau}
   (w^{- 1})_{\tau} . \]
We prove that for $L^2$ integrable Hilbert space valued martingales $Y, X$
with $Y$ differentially subordinate to $X$ there holds
\[ \| Y \|_{L^2 (w)} \lesssim Q^{\mathcal{F}}_2 [w] \| X \|_{L^2 (w)} \]
where the implied constant is numeric and does not depend upon the dimension,
the pair of martingales or the weight. The linear growth in the quantity $Q_2
(w)$ is sharp.

\

The proof in this paper is different from the proofs in {\cite{Lac2015a}} and
{\cite{ThiTreVol2015a}}. In {\cite{Lac2015a}} so--called sparse operators are
used while in {\cite{ThiTreVol2015a}} the authors reduce the estimate through
the use of so--called outer measure space theory.

Our approach is the following. We derive an explicit Bellman function of four
variables adapted to the problem. It has certain conditions on its range, a
continuous sub--convexity as well as discrete one--leg convexity, such as seen
in {\cite{ThiTreVol2015a}} for two smaller Bellman functions (their functions
make up a part of ours). We heavily use the explicit form of our Bellman
function and its regularity properties in several parts in our proof to handle
the delicacy of the continuous--in--time processes with values in Hilbert
space. The resulting function is in the `dualized' or `weak form', which is in
a contrast to the `strong form' of a Burkholder type functional often seen
when using the strong subordination condition (\ref{defsub}). (The explicit
form of a Burkholder type functional for this weighted question is still
open). Indeed, the form of the strong differential subordination condition is
adapted to work well for Burkholder type functionals and arises naturally in
this setting. The passage to its use in the weak form is accomplished through
the use of the so--called Ellipse Lemma and requires a Bellman function
solving the entire problem at once as opposed to splitting the problem into
pieces. This is the first use of this strategy for problems in probability and
should allow generalisations of numerous existing results as well as an
alternative (allbeit more complicated) proof of Wang's extension to
Burkholder's famous estimates using {\cite{VasVol2010c}} or
{\cite{BanOse2014b}}. Note that for these $L^p$ problems, fewer difficulties
arise, even in the presence of jumps. This is thanks to the convexity of the
domain in the $L^p$ problem. The discrete convexity required to control the
jumps is almost free, when using a trick from {\cite{DomPet2014c}}. This trick
is not available here because of the non--convex domain.

\

Our result gives through the formula in {\cite{BanJan2008a}} a probabilistic
proof of the weighted estimate for the Beurling--Ahlfors transform with its
implication, a famous borderline regularity problem for the Beltrami equation,
solved in {\cite{PetVol2002a}}. Other applications are discussed in the last
section. They include a dimensionless weighted bound for discrete and
semi-discrete second order Riesz transforms.

\

\subsection{Differentially subordinate martingales}

\

Consider first discrete--in--time martingales. For that let $(\Omega,
\mathcal{F}_{\infty}, \mathbbm{P})$ a probability space with a nondecreasing
sequence $\mathcal{F} = (\mathcal{F}_n)_{n \geqslant 0}$ of sub
$\sigma$--fields of $\mathcal{F}_{\infty}$ such that $\mathcal{F}_0$ contains
all $\mathcal{F}_{\infty}$--null sets. We are interested in
$\mathbbm{H}$--valued martingales, where $\mathbbm{H}$ is a separable Hilbert
space with norm $\left| \!\! \; \cdot \; \right|_{\mathbbm{H}}$ and scalar
product $\langle \cdot, \cdot \rangle_{\mathbbm{H}}$: if $f = \{ f_n \}_{n \in
\mathbbm{N}}$ is a $\mathbbm{H}$--valued martingale adapted to $\mathcal{F}$,
we note $f_n = \sum_{k = 0}^n \mathd f_k$, with the convention $\mathd f_0
\assign f_0$, and $\mathd f_k \assign f_k - f_{k - 1}$, for $k \geqslant 1$.
Similarily, if $g$ is another adapted $\mathbbm{H}$--valued martingale, we
note $g_n = \sum_{k = 0}^n \mathd g_k$ with the same conventions. One says
that $g$ is differentially subordinate to $f$ if one has for almost all
$\omega \in \Omega$ and all $k \geqslant 0$, $| \mathd g_k |_{\mathbbm{H}}
\leqslant | \mathd f_k |_{\mathbbm{H}}$.

\

In this paper we consider continuous--in--time filtrations. Let again
$(\Omega, \mathcal{F}_{\infty}, \mathbbm{P})$ a probability space with a
nondecreasing right continuous family $\mathcal{F}= (\mathcal{F}_t)_{t
\geqslant 0}$ of sub $\sigma$--fields of $\mathcal{F}_{\infty}$ such that
$\mathcal{F}_0$ contains all $\mathcal{F}_{\infty}$--null sets. We are
interested in $\mathbbm{H}$--valued c{\`a}dl{\`a}g martingales, where
$\mathbbm{H}$ is a separable Hilbert space. In order to clearly define
differential subordination in this setting, we make use of the square bracket
or quadratic variation process.

\

Recall that the quadratic variation process of a semimartingale $X$ is the
process denoted by $[X, X] \assign ([X, X]_t)_{t \geqslant 0}$ and defined as
(see e.g. {\tmname{Protter}} {\cite{Pro2005a}})
\[ [X, X]_t = X_t^2 - 2 \int_0^t X_{s -} \mathd X_s \]
where we have set $X_{0 -} = 0$. Similarily, the quadratic covariation of two
semimartingales $X$ and $Y$ is the following process also known as the bracket
process
\[ [X, Y]_t \assign X_t Y_t - \int^t_0 X_{s -} \mathd Y_s - \int Y_{s -}
   \mathd X_s . \]

\begin{definition}[differential subordination]
Let $X$ and $Y$ two adapted c{\`a}dl{\`a}g semimartingales
taking values in a separable Hilbert space. We say $Y$ is differentially
subordinate by quadratic variation to $X$ iff $$[X, X]_t - [Y, Y]_t$$ is a
nondecreasing and nonnegative function of $t \geqslant
0$.
\end{definition}

Let us denote by $X^c$ the unique continuous part of $X$ with
\[ [X, X]_t = | X_0 |^2 + [X^c, X^c]_t + \sum_{0 < s \leqslant t} | \Delta X_s
   |^2 . \]
There holds $[X, X]^c_t = [X^c, X^c]_t$ and $\Delta [X, X]_t = | \Delta X_t
|^2$ where $\Delta X_t \assign X_t - X_{t -}$. We have the following obvious
characterisation distinguishing the continuous and jump parts:

\begin{lemma}
If $X$ and $Y$ are semimartingales,
then $Y$ is differentially subordinate to $X$ if and only if (i) $[X, X]^c_t -
[Y, Y]^c_t$ is a nonnegative and nondecreasing function of $t$, (ii) the
inequality $| \Delta Y_t | \leqslant | \Delta X_t |$ holds for all $t > 0$ and
(iii) $| Y_0 | \leqslant | X_0 |$.
\end{lemma}

\

\subsection{Martingales in non--homogeneous weighted spaces}

\

Let again $(\Omega, \mathcal{F}_{\infty}, \mathbbm{P})$ a probability space
with a nondecreasing right continuous family $\mathcal{F} \assign
(\mathcal{F}_t)_{t \geqslant 0}$ of sub $\sigma$--fields of
$\mathcal{F}_{\infty}$ such that $\mathcal{F}_0$ contains all
$\mathcal{F}_{\infty}$--null sets. The measure $\mathd \mathbbm{P}$ is
arbitrary (up to the obvious normalisation). If $X$ and $Y$ are adapted
c{\`a}dl{\`a}g square integrable $\mathbbm{H}$--valued martingales and $Y$ is
differentially subordinate to $X$, then it is obvious that
\begin{equation}
  \| Y \|_2 \leqslant \| X \|_2 . \label{eq: L2 estimates without change of
  law}
\end{equation}
Recall here that $\| X \|_2 \assign \sup_t \| X_t \|_2$, where
\begin{equation}
  \| X_t \|^2_2 \assign \mathbbm{E} | X_t |^2 = \int_{\Omega} | X_t (\omega)
  |^2 \mathd \mathbbm{P} (\omega) . \label{eq: definition L2 norm of
  semimartingale}
\end{equation}

Assume again that $Y$ is differentially subordinate to $X$. We might insist on
the underlying probability space at hand by saying in short that $X$ and $Y$
are $\mathbbm{P}$--martingales and that $Y$ is
$\mathbbm{P}$--differentially--subordinate to $X$. The main concern of this
paper is to obtain sharp inequalities similar to (\ref{eq: L2 estimates
without change of law}) under a change of law in the definition of the
$L^2$--norm according to {\cite{DelMey1982a}}. Let $w$ be a positive, uniformy
integrable martingale (that we often identify with its closure $w_{\infty}$)
that we call a weight. Let $\mathd \mathbbm{Q} \assign \mathd \mathbbm{P}^w
\assign w \mathd \mathbbm{P}$ and $(\Omega, \mathcal{F}_{\infty},
\mathbbm{Q})$ be a probability space with the same assumptions as $(\Omega,
\mathcal{F}_{\infty}, \mathbbm{P})$ but with a change of the probability law.

\begin{question}
Let $\mathbbm{P}$ and
$\mathbbm{Q}$ such that $(\Omega, \mathcal{F}_{\infty}, \mathbbm{P})$ and
$(\Omega, \mathcal{F}_{\infty}, \mathbbm{Q})$ are two filtered probability
spaces as described above. \ Does there exist a constant $C_{\mathbbm{P},
\mathbbm{Q}} > 0$ depending only on $(\mathbbm{P}, \mathbbm{Q})$ such that if
$X$ and $Y$ are uniformly integrable $\mathbbm{P}$--martingales adapted to
$\mathcal{F}$ and $Y$ is $\mathbbm{P}$--differentially--subordinate to $X$,
then
\[ \| Y \|_{2, \mathbbm{Q}} \leqslant C_{\mathbbm{P}, \mathbbm{Q}} \| X \|_{2,
   \mathbbm{Q}} . \]
   \end{question}

We look for $C_{\mathbbm{P}, \mathbbm{Q}} \assign C_w$ allowing to compare $\|
Y \|_{2, \mathbbm{Q}} \assign \| Y \|_{2, w}$ and $\| X \|_{2, \mathbbm{Q}}
\assign \| X \|_{2, w}$. We will need also $u = w^{- 1}$ the inverse weight
and we assume $u$ uniformly integrable. We will finally note $\mathd
\mathbbm{P}^u \assign u \mathd \mathbbm{P}$. It follows that $\mathbbm{P}^w$
and $\mathbbm{P}^u$ are probability measures on $\Omega$ up to the obvious
normalisations. The necessary condition on the weight is classical:

\begin{definition}[$A_2$ class]
Let
$(\Omega, \mathcal{F}, (\mathcal{F}_t)_{t \geqslant 0}, \mathbbm{P})$ a
filtered probability space. We say that the weight $w > 0$ is in the
$A_2$--class, iff the $A_2$--characteristic of the weight $w$, noted
$Q^{\mathcal{F}}_2 [w]$ and defined as
\[ Q^{\mathcal{F}}_2 [w] \assign \sup_{\tau}  \tmop{ess} .
   \sup_{\omega}  (w_{})_{\tau} (w^{- 1})_{\tau} \]
with the first supremum running over all adapted stopping times, is
finite.
\end{definition}

We often write $Q^{\mathcal{F}}_2 [w] \assign \sup_{\tau} \tmop{ess} .
\sup_{\omega} w_{\tau} u_{\tau}$ where $u \assign w^{- 1}$ is the
inverse weight.

\

\section{Statement of the main results}

\

\begin{theorem}[differential subordination under change of law]
  \label{T: main result}Let $X$ and $Y$ be two adapted uniformly integrable
  c{\`a}dl{\`a}g $\mathbbm{H}$--valued martingales such that $Y$ is
  differentially subordinate to $X$. Let $w$ be an admissible weight in the
  $\tmmathbf{A}_2$ class. Then
  \[ \| Y \|_{L^2 (w)} \lesssim Q^{\mathcal{F}}_2 [w]  \| X \|_{L^2 (w)} \]
  and the linear growth in $Q^{\mathcal{F}}_2 [w]$ is sharp.
\end{theorem}

\

This result will be a consequence of the following bilinear estimate:

\begin{proposition}[bilinear estimate]
  \label{P: bilinear estimate}Let $X$ and $Y$ be two adapted uniformly
  integrable c{\`a}dl{\`a}g $\mathbbm{H}$--valued martingales such that $Y$ is
  differentially subordinate to $X$. Let $w$ an admissible weight in the
  $\tmmathbf{A}_2$ class. Then
  \[ \mathbbm{E} \int_0^{\infty} | \mathd [Y, Z]_t | \lesssim
     Q^{\mathcal{F}}_2 [w]  \| X \|_w \| Z \|_u . \]
\end{proposition}

We have an explicit expression of the function described below. This is, aside
from Theorem \ref{T: main result}, one of the main results of this paper. Let
us note $V$ the quadruplet
\[ V \assign (x, y, r, s) \in \mathbbm{H} \times \mathbbm{H} \times
   \mathbbm{R}^{\ast}_+ \times \mathbbm{R}^{\ast}_+ = : \mathbbm{S}. \]
The variables $(x, y)$ will be associated to $\mathbbm{H}$--valued martingales
whereas the variables $(r, s)$ to $\mathbbm{R}$--valued martingales for the
weights. We introduce $\mathcal{D}_Q$ the domain
\[ \mathcal{D}_Q \assign \left\{ V \in \mathbbm{S}: \hspace{1em} 1 \leqslant r
   s \leqslant Q \right\} . \]
We will often restrict our attention to truncated weights, that is given $0 <
\varepsilon < 1$, variables $(r, s)$ bounded below and above
\[ \mathcal{D}^{\varepsilon}_Q \assign \left\{ V \in \mathcal{D}_Q :
   \hspace{1em} \varepsilon \leqslant r \leqslant \varepsilon^{- 1},
   \hspace{1em} \varepsilon \leqslant s \leqslant \varepsilon^{- 1} \right\} .
\]
\begin{lemma}[existence and properties of the Bellman function]
  \label{L: existence and properties of the Bellman function}There exists a
  function $B (V) = B_Q$ that is $\mathcal{C}^1$ on
  $\mathcal{D}^{\varepsilon}_Q$, and piecewise $\mathcal{C}^2$, with the
  estimate
  \[ B (V) \lesssim \frac{| x |^2}{r} + \frac{| y |^2}{s} \]
  and on each subdomain where it is $\mathcal{C}^2$ there holds
  \[ \mathd^2 B \geqslant \frac{2}{Q} | \mathd x | | \mathd y | . {\color{red}
     } \]
  Whenever $V$ and $V_0$ are in the domain, the function has the property
  \[ B (V) - B (V_0) - \mathd B (V_0) (V - V_0) \geqslant \frac{2}{Q} | x -
     x_0 |  | y - y_0 | . \]
  Moreover, we have the estimates
  \[ | (\partial^2_x B \mathd x, \mathd x) | \lesssim \varepsilon^{- 1} |
     \mathd x |^2, \hspace{2em} | (\partial^2_y B \nocomma \mathd y, \mathd y)
     | \lesssim \varepsilon^{- 1} | \mathd y |^2 \]
  with the implied constants independent of $V$ and $(\mathd x, \mathd y)$.
\end{lemma}

\

\section{Existence and properties of the Bellman function}

\

\begin{proof*}{Proof of Lemma \ref{L: existence and properties of the Bellman
function}}{\dueto{existence and properties of the Bellman function}}
  
  We give an
  explicit expression for such a function. Let $V = (x, y, r, s)$ and $W = (r,
  s)$. We first consider

  \[ B_1 (x, y, r, s) = \frac{\langle x, x \rangle}{r} + \frac{\langle y, y
     \rangle}{s} . \]
  Then trivially $0 \leqslant B_1 \leqslant \frac{\langle x, x \rangle}{r} +
  \frac{\langle y, y \rangle}{s}$ and
  \begin{eqnarray*}
    (\mathd^2 B_1 \mathd V \nocomma, \mathd V) & = & \frac{2}{r} \langle
    \mathd x, \mathd x \rangle + \frac{2 \langle x, x \rangle}{r^3} (\mathd
    r)^2 - 4 \frac{\langle x, \mathd x \rangle}{r^2} \mathd r\\
    &  &\hspace{1em} + \frac{2}{s} \langle \mathd y, \mathd y \rangle + \frac{2 \langle
    y, y \rangle}{s^3} (\mathd s)^2 - 4 \frac{\langle y, \mathd y
    \rangle}{s^2} \mathd s\\
    & = & \frac{2}{r} \left\langle \mathd x - \frac{x}{r} \mathd r, \mathd x
    - \frac{x}{r} \mathd r \right\rangle \\
    &&\hspace{1em}+ \frac{2}{s} \left\langle \mathd y -
    \frac{y}{s} \mathd s, \mathd y - \frac{y}{s} \mathd s \right\rangle\\
    & \geqslant & 0
  \end{eqnarray*}
  Letting $V_0 = (x_0, y_0, r_0, s_0)$ and $V = (x, y, r, s)$ also calculate
  \begin{eqnarray*}
    &  & - (B_1 (V_0) - B_1 (V) + \mathd B_1 (V_0) (V - V_0))\\
    & = & - \left( \frac{x^2_0}{r_0} - \frac{x^2}{r} + \frac{2 x_0}{r_0} (x -
    x_0) - \frac{x^2_0}{r^2_0} (r - r_0) \right)\\
    &  & \hspace{1em}- \left( \frac{y^2_0}{s_0} - \frac{y^2}{s} + \frac{2 y_0}{s_0} (y -
    y_0) - \frac{y^2_0}{s^2_0} (s - s_0) \right)\\
    & = & r \left\langle \frac{x}{r} - \frac{x_0}{r_0}, \frac{x}{r} -
    \frac{x_0}{r_0} \right\rangle + s \left\langle \frac{y}{s} -
    \frac{y_0}{s_0}, \frac{y}{s} - \frac{y_0}{s_0} \right\rangle .
  \end{eqnarray*}
  We now consider the two functions from {\cite{ThiTreVol2015a}}
  \[ K (r, s) = \frac{\sqrt{r s}}{\sqrt{Q}} \left( 1 - \frac{\sqrt{r s}}{8
     \sqrt{Q}} \right) \]
   \[  N (r, s) = \frac{\sqrt{r s}}{\sqrt{Q}}
     \left( 1 - \frac{(r s)^2}{128 Q^2} \right) \]
  in the domain $1 \leqslant r s \leqslant Q$. We have
  \[ 0 \leqslant K \leqslant \left( 1 - \frac{1}{8 \sqrt{Q}} \right)
     \frac{\sqrt{r s}}{\sqrt{Q}} < \frac{\sqrt{r s}}{\sqrt{Q}} \leqslant 1, \]
  \[ 0 \leqslant N \leqslant \left( 1 - \frac{1}{128 Q^2} \right)
     \frac{\sqrt{r s}}{\sqrt{Q}} < \frac{\sqrt{r s}}{\sqrt{Q}} \leqslant 1. \]
  So in particular $r s \geqslant r s - K^2 > r s \left( 1 - \frac{1}{Q}
  \right)$. One calculates that
  \[ - (\mathd^2 K \mathd W, \mathd W) \geqslant \frac{1}{8 Q} | \mathd r | |
     \mathd s |, \]
  \[ - (\mathd^2 N \mathd W, \mathd W) \gtrsim \frac{1}{Q^2} s^2 (\mathd r)^2, \]
  \[ - (\mathd^2 N \mathd W, \mathd W) \gtrsim \frac{1}{Q^2} r^2
     (\mathd s)^2 . \]
  \[ \  \]
  One also has whenever $W, W_0$ in the domain then
  \[ K (W_0) - K (W) + \mathd K (W_0) (W - W_0) \gtrsim \frac{1}{Q} | r - r_0
     |  | s - s_0 |, \]
  \[ N (W_0) - N (W) + \mathd N (W_0) (W - W_0) \gtrsim \frac{1}{Q^2} s_0 s |
     r - r_0 |^2, \]
  \[ N (W_0) - N (W) + \mathd N (W_0) (W - W_0) \gtrsim \frac{1}{Q^2} r_0 r |
     s - s_0 |^2 . \]
  These remarkable one-leg concavity properties were proven in {\cite{ThiTreVol2015a}}.
  
  \
  
  Let now
  \[ B_2 = \frac{\langle x, x \rangle}{2 r - \frac{1}{s (N (r, s) + 1)}} +
     \frac{\langle y, y \rangle}{s}= \frac{\langle x, x \rangle}{r + M (r, s)} +
  \frac{\langle y, y \rangle}{s} , \]
   where 
   \[M (r, s) = r - \frac{1}{s (N (r, s) + 1)} . \] 
   One checks easily by calculation of their Hessians that
  \[ F (x, r, M) = \frac{\langle x, x \rangle}{r + M} ,\]
  \[G (r, s, N) =
     \frac{1}{s (N + 1)} \]
  are convex everywhere. In order to estimate the Hessian of $B_2$ from below,
  one merely requires estimates of derivatives
  \[ - \partial_M F = \frac{\langle x, x \rangle}{(r + M)^2} \geqslant
     \frac{\langle x, x \rangle}{4 r^2} \tmop{and} - \partial_N G = \frac{1}{s
     (N + 1)^2} \geqslant \frac{1}{4 s} . \]
  Since $0 \leqslant r - \frac{1}{s (N (r, s) + 1)} \leqslant r$ we know $0
  \leqslant B_2 \leqslant \frac{| x |^2}{r} + \frac{| y |^2}{s}$. Now the
  Hessian estimate becomes
  \begin{eqnarray*}
    &&(\mathd^2 B_2 \mathd V, \mathd V)\\
     & \gtrsim & \frac{\langle x, x
    \rangle}{4 r^2} \frac{1}{s (N + 1)^2} \frac{1}{Q^2} | \mathd r |^2 s^2 +
    \frac{2}{s} \left\langle \mathd y - \frac{y}{s} \mathd s, \mathd y -
    \frac{y}{s} \mathd s \right\rangle\\
    & \gtrsim & \frac{| x |^2 s}{Q^2 r^2} | \mathd r |^2 + \frac{2}{s}
    \left\langle \mathd y - \frac{y}{s} \mathd s, \mathd y - \frac{y}{s}
    \mathd s \right\rangle\\
    & \gtrsim & \frac{| x |}{Q} | \mathd r |  \left| \mathd y - \frac{y}{s}
    \mathd s \right|.
  \end{eqnarray*}
  This function has the additional property
  \begin{eqnarray*}
  && - (B_2 (V_0) - B_2 (V) + \mathd B_2 (V_0) (V - V_0)) \\
  &\gtrsim &
     \frac{\langle x_0, x_0 \rangle}{Q^2 r^2_0} s (r - r_0)^2 + s \left\langle
     \frac{y}{s} - \frac{y_0}{s_0}, \frac{y}{s} - \frac{y_0}{s_0}
     \right\rangle .
  \end{eqnarray*}   
  Indeed, write
  \[ \frac{\langle x, x \rangle}{2 r - \frac{1}{s (N (r, s) + 1)}} = H (x, r,
     s, N (r, s)) \tmop{with} H (x, r, s, N) = \frac{\langle x, x \rangle}{2 r
     - \frac{1}{s (N + 1)}} \]
  where $H$ is convex and
  \[ - \partial_N H \gtrsim \frac{\langle x, x \rangle}{Q^2 r^2 s} . \]
  Now since $H$ is convex we have with $P_0 = (x_0, r_0, s_0, N_0)$ and with $P =
  (x, r, s, N)$ that $H (P) \geqslant H (P_0) + \mathd H (P - P_0)$. So
  \begin{eqnarray*}
    &  & H (P) - H (P_0)  \\
    &  & \hspace{1em} - \partial_x H (P_0) (x - x_0) - \partial_r H (P_0)
    (r - r_0) - \partial_s H (P_0) (s - s_0)\\
    & \geqslant & - \partial_N H (P_0) (N_0 - N) .
  \end{eqnarray*}
  With
  \begin{eqnarray*}
    &  & N (r_0, s_0) - N (r, s) + \partial_r N (r_0, s_0) (r - r_0)  + \partial_s
     N (r_0, s_0) (s - s_0) \\
     &\gtrsim & \frac{1}{Q^2} r_0 r | s - s_0 |^2 
  \end{eqnarray*}
  the above becomes with $N_0 = N (r_0, s_0)$ and $N = N (r, s)$
  \begin{eqnarray*} 
  && B_2 (V) - B_2 (V_0) - \mathd B_2 (V_0) (V - V_0)\\
  & \gtrsim & \frac{\langle
     x_0, x_0 \rangle}{Q^2 r^2_0} s | r - r_0 |^2 + s \left\langle \frac{y}{s}
     - \frac{y_0}{s_0}, \frac{y}{s} - \frac{y_0}{s_0} \right\rangle 
   \end{eqnarray*}
  where we used the lower derivative estimate and the chain rule. Analogously
  \[ B_3 = \frac{\langle x, x \rangle}{r} + \frac{\langle y, y \rangle}{2 s -
     \frac{1}{r (N (r, s) + 1)}} \]
  has the same size estimates as well as
  \[ (\mathd^2 B_3 \mathd V, \mathd V) \gtrsim \frac{| y |}{Q} | \mathd s | 
     \left| \mathd x - \frac{x}{r} \mathd r \right| \]
  and one--leg convexity
  \begin{eqnarray*}
  && B_3 (V) - B_3 (V_0) - \mathd B_3 (V_0) (V - V_0) \\
  &\gtrsim &\frac{\langle
     y_0, y_0 \rangle}{Q^2 s^2_0} r | s - s_0 |^2 + r \left\langle \frac{x}{r}
     - \frac{x_0}{r_0}, \frac{x}{r} - \frac{x_0}{r_0} \right\rangle . 
   \end{eqnarray*}

  Let us now consider
  \[ H_4 (x, y, r, s, K) = \sup_{0 < a} \beta (a, x, y, r, s, K) = \sup_{0 <
     a} \left( \frac{\langle x, x \rangle}{r + a K} + \frac{\langle y, y
     \rangle}{s + a^{- 1} K} \right) . \]
  Testing for critical points gives 
  \[\partial_a \beta = - \frac{\langle x, x
  \rangle K}{(r + a K)^2} + \frac{\langle y, y \rangle K}{(a s + K)^2}.\] 
  So
  $\partial_a \beta = 0$ if and only if 
  \[ a = a' = \frac{| y | r - | x | K}{| x
  | s - | y | K}.\] 
  Since only $a > 0$ are admissible, we require that $| y | r
  - | x | K$ and $| x | s - | y | K$ have the same sign. To determine sign
  change of $\partial_a \beta$ at $a'$, Consider
  \[ - \frac{| x |}{r + a K} + \frac{| y |}{a s + K} = \frac{(| y | r - | x |
     K) - a (| x | s - | y | K)}{(r + a K) (a s + K)} . \]
  If the signs are negative, then the sign change is from negative to positive
  otherwise from positive the negative. For a maximum to be attained at $a' >
  0$ we require that both numerator and denominator be positive. Then, if $K$
  is relatively small, meaning $| y | r - | x | K$ and $| x | s - | y | K$
  positive we have
  \begin{eqnarray*}
    H_4 (x, y, r, s, K) & = & \beta (a', x, y, r, s)\\
    & = & \frac{\langle x, x \rangle (| x | s - | y | K)}{r (| x | s - | y |
    K) + (| y | r - | x | K) K}\\
    &  & \hspace{1em}+ \frac{\langle y, y \rangle (| y | r - | x | K)}{s (| y | r - | x |
    K) + (| x | s - | y | K) K}\\
    & = & \frac{\langle x, x \rangle s - 2 | x | | y | K + \langle y, y
    \rangle r}{r s - K^2} .
  \end{eqnarray*}
  Observe that by the above considerations on $K$, the denominator is never 0.
  The case $| x | = 0$ or $| y | = 0$ corresponds to other parts of the
  domain, so when $K$ is small in the sense above, this function is in
  $\mathcal{C}^2$.
  
  When $| y | r - | x | K \leqslant 0$ or $| x | s - | y | K \leqslant 0$, the
  supremum is attained at the boundary and $H_4 = \frac{\langle y, y
  \rangle}{s}$ or $H_4 = \frac{\langle x, x \rangle}{r}$. Thanks to the size
  restrictions on $K$ we never have both $| x | s - | y | K \leqslant 0$ and
  $| y | r - | x | K \leqslant 0$ unless $x, y = 0$, indeed
  \begin{eqnarray*}
   && | x | (| x | s - | y | K) + | y | (| y | r - | x | K) \\
   &=& \frac{| x
    |^2}{r} - 2 \frac{| x | | y |}{r s} K + \frac{| y |^2}{s}\\
    & = & \left( \frac{| x |}{\sqrt{r}} - \frac{| y |}{\sqrt{s}} \right)^2 +
    \frac{2 | x | | y |}{\sqrt{r s}} \left( 1 - \frac{K}{\sqrt{r s}} \right).
  \end{eqnarray*}
  With $1 - \frac{K}{\sqrt{r s}} > 0$ we see that the above is never negative
  and the quantity vanishing implies $x = y = 0$. If $| x | s - | y | K
  \leqslant 0$ and $| y | r - | x | K > 0$ then $\frac{\langle x, x
  \rangle}{r} < \frac{\langle y, y \rangle}{s}$ and $H_4 = \frac{\langle y, y
  \rangle}{s}$, if $| y | r - | x | K \leqslant 0$ and \ $| x | s - | y | K >
  0$ then $H_4 = \frac{\langle x, x \rangle}{r}$.
  
  Notice that when $\frac{\langle x, x \rangle}{r}$=$\frac{\langle y, y
  \rangle}{s}$ and $x, y \neq 0$ then $| y | r - | x | K > 0$ and $| x | s - |
  y | K > 0$. Indeed, we have seen $\frac{| x |^2}{r} - 2 \frac{| x | | y |}{r
  s} K + \frac{| y |^2}{s} > 0$. Thus $\frac{\langle x, x \rangle}{r} =
  \frac{\langle y, y \rangle}{s} > \frac{| x | | y |}{r s} K$ and $| y | r - |
  x | K > 0$ and $| x | s - | y | K > 0$.
  
  Thus $H_4 \in \mathcal{C}^2$ for these parts of the domain. We also see from
  these considerations that in order to see $H_4 \in \mathcal{C}^1$ we only
  need to check the cuts $| x | s - | y | K = 0$ and $| y | r - | x | K
  \geqslant 0$ as well as $| y | r - | x | K = 0$ and $| x | s - | y | K
  \geqslant 0$.
  
  When $| y | r - | x | K > 0$ and $| x | s - | y | K > 0$ (we call this part
  of the domain $R_1$)
  \[ (\partial_x H_4, \mathd x) = 2 \frac{\langle d x, x \rangle}{| x |}
     \frac{| x | s - | y | K}{r s - K^2} \]
  \[ \partial_r H_4 = - \frac{(| x | s - | y | K)^2}{(r s - K^2)^2} \]
  \[ \nocomma \partial_K H_4 = - 2 \frac{(| x | s - | y | K) (| y | r - | x |
     K)}{(r s - K^2)^2} \]
  We first prove that $\partial_x H_4$ is continuous throughout. Recall that
  we have to treat three regions: $R_1$ and $R_2$ where $| y | r - | x | K >
  0$ and $| x | s - | y | K \leqslant 0$ and $R_3$ where $| x | s - | y | K
  \leqslant 0$ and $| y | r - | x | K > 0$. Inside $R_2$ we have $H_4 =
  \frac{\langle y, y \rangle}{s}$ and thus $\partial_x H_4 = 0$. Inside $R_3$
  we have $H_4 = \frac{\langle x, x \rangle}{r}$ and thus $\partial_x H_4 =
  \frac{2 \langle x, d x \rangle}{r}$. Inside $R_1$
  \begin{eqnarray*}
 \partial_x H_4 &=& 2 \frac{\langle \mathd x, x \rangle}{| x |} \frac{| x |
     s - | y | K}{r s - K^2} \\
     &=& 2 \langle x, \mathd x \rangle \frac{| x | s - |
     y | K}{r (| x | s - | y | K) + (| y | r - | x | K) K} . \end{eqnarray*}
  We have three cases, first, let us approach a boundary point of $R_1$ from
  within $R_1$ so that $| y | r - | x | K > 0$ and $| x | s - | y | K = 0$.
  Assume therefore $| y | r - | x | K \sim a > 0$ and $0 < | x | s - | y | K <
  \varepsilon$. There holds $| \langle \partial_x H_4, \mathd x \rangle |
  \leqslant 2 | \mathd x | \frac{\varepsilon}{r s - K^2} \lesssim \varepsilon
  | \mathd x |$ since $r s - K^2$ is bounded below. Letting $\varepsilon
  \rightarrow 0$ shows continuity in this point. Second, let us approach a
  boundary point \ $| x | s - | y | K > 0$ and $| y | r - | x | K = 0$ from
  within $R_1$. Assume therefore $| x | s - | y | K \sim a > 0$ and $0 < | y |
  r - | x | K < \varepsilon$. We show there holds $(\partial_x H_4, d x)
  \lesssim \frac{\varepsilon}{a} | \mathd x |$. Since
  \begin{eqnarray*}
  && \frac{1}{r} - (| y | r - | x | K) K \frac{| x | s - | y | K}{r^2 (| x | s
     - | y | K)^2} \\
     &\leqslant & \frac{(| x | s - | y | K)}{r (| x | s - | y | K)
     + (| y | r - | x | K) K}\\
     & \leqslant &\frac{1}{r} 
 \end{eqnarray*}
  we have
  \begin{eqnarray*}
   && \left| \frac{2 \langle x, \mathd x \rangle (| x | s - | y | K)}{r (| x | s
    - | y | K) + (| y | r - | x | K) K} - \frac{2 \langle x, \mathd x
    \rangle}{r} \right| \\
    & \leqslant & 2 | \langle x, \mathd x \rangle |
    \frac{(| y | r - | x | K) K}{r^2 (| x | s - | y | K)}\\
    & \lesssim & | x | | \mathd x | \frac{\varepsilon}{a} .
  \end{eqnarray*}
  Since $0 < | y | r - | x | K < \varepsilon$ and $s, r, K$ controlled, one
  can deduce from $| x | s - | y | K \sim a$ that $| x | \sim a$. Last, let us
  approach $| y | r - | x | K = 0$ and $| x | s - | y | K = 0$. To this end,
  one can see that if $0 < | y | r - | x | K < \varepsilon$ and $0 < | x | s -
  | y | K < \varepsilon$ then $| x |, | y | \lesssim \varepsilon$,
  establishing continuity in the third case.
  
  The $\partial_r H_4$ derivative is similar since the term $\frac{| x | s - |
  y | K}{r s - K^2}$ reappears as a square and in $R_3$ notice that $H_4 =
  \frac{\langle x, x \rangle}{r}$ so $\partial_r H_4 = - \frac{\langle x, x
  \rangle}{r^2}$. It is easy to see that the derivative $\partial_K H_4$ is
  zero in $R_2$ and $R_3$ as well as when approaching the boundary of $R_1$.
  
  These derivatives are representative by symmetry and the function is
  therefore in $\mathcal{C}^1$. As a consequence
  \[ B_4 (X, Y, x, y, r, s) = H_4 (x, y, r, s, K (r, s)) \in \mathcal{C}^1 .
  \]
  Function $B_4$ is as supremum of convex functions convex. It has been shown
  indirectly in {\cite{NazTreVol1999a}} that $- \partial_K B_4 \geqslant 0$
  everywhere and that in $R'_1 \subset R_1$ where $| y | r - 2 | x | K > 0$
  and $| x | s - 2 | y | K > 0$ we have $- \partial_K B_4 \gtrsim \frac{| x |
  | y |}{r s}$. We present an easier argument. Recall that
  \begin{eqnarray*}
    &&- \partial_K B_4 \\
    & = & 2 \frac{(| x | s - | y | K) (| y | r - | x | K)}{(r
    s - K^2)^2}\\
    & = & 2 \frac{(| x | s - | y | K) (| y | r - | x | K) | x | | y |}{(r (|
    x | s - | y | K) + K (| y | r - | x | K)) (s (| y | r - | x | K) + K (| x
    | s - | y | K))}
  \end{eqnarray*}
  So $- \partial_K B_4 \geqslant c \frac{| x | | y |}{r s}$ if
  \[ \frac{r s}{c} \geqslant K^2 \noplus + r s + \frac{K r (| x | s - | y |
     K)}{| y | r - | x | K} + \frac{K s (| y | r - | x | K)}{| x | s - | y |
     K} .\]
  Now $K^2 \leqslant 1 \leqslant r s$ and when $| y | r - 2 | x | K \geqslant
  0$ then $| y | r - | x | K \geqslant | x | K$. Similarly $| x | s - | y | K
  \geqslant | y | K$. So the last two terms are bounded by $\frac{K r | x |
  s}{| x | K} + \frac{K s | y | r}{| y | K} = 2 r s$. So $c = 1 / 4$ works. In
  $R'_1$
  \[ (\mathd^2 B_4 \mathd V, \mathd V) \geqslant 4 \frac{| x | | y |}{8 r s Q}
     | \mathd r | | \mathd s | = \frac{| x | | y |}{2 r s Q} | \mathd r | |
     \mathd s | . \]
  We need to add more functions with the good concavity for other $K$. Let
  \[ B_5 = \frac{\langle x, x \rangle}{2 r - \frac{1}{s (K (r, s) + 1)}} +
     \frac{\langle y, y \rangle}{s} . \]
  Since $0 \leqslant r - \frac{1}{s (K (r, s) + 1)} \leqslant r$ we know $0
  \leqslant B_5 \leqslant \frac{| x |^2}{r} + \frac{| y |^2}{s}$. Now the
  Hessian estimate becomes
  \[ (\mathd^2 B_5 \mathd V, \mathd V) \geqslant \frac{\langle x, x \rangle}{4
     r^2} \frac{1}{s (K + 1)^2} \frac{1}{8 Q} | \mathd r | | \mathd s |
     \geqslant \frac{| x |^2}{128 Q s r^2} | \mathd r | | \mathd s | . \]
  $B_5$ convex and when $2 | x | K \geqslant | y | r$ then
  \[ (\mathd^2 B_5 \mathd V, \mathd V) \geqslant \frac{| x | | y |}{256 K Q s
     r} | \mathd r | | \mathd s | \geqslant \frac{| x | | y |}{256 Q s r} |
     \mathd r | | \mathd s | . \]
  With
  \[ B_6 = \frac{\langle x, x \rangle}{r} + \frac{\langle y, y \rangle}{2 s -
     \frac{1}{r (K (r, s) + 1)}} \]
  we have $0 \leqslant B_6 \leqslant \frac{| x |^2}{r} + \frac{| y |^2}{s}$
  convex and when $2 | y | K \geqslant | x | s$ then
  \[ (\mathd^2 B_6 \mathd V, \mathd V) \geqslant \frac{| x | | y |}{256 Q s r}
     | \mathd r | | \mathd s | . \]
  Together, we have for $B_7 = B_4 + B_5 + B_6$ that
  \[ (\mathd^2 B_7 \mathd V, \mathd V) \gtrsim \frac{| x | | y |}{Q s r} |
     \mathd r | | \mathd s | . \]
  Through similar considerations as above, we have discrete one-leg convexity
  \[ B_7 (V) - B_7 (V_0) - \mathd B_7 (V_0) (V - V_0) \gtrsim \frac{| x_0 | |
     y_0 |}{Q s_0 r_0} | r - r_0 | | s - s_0 | . \]

  Letting for appropriate fixed $c_i$
  \begin{equation}
    B = c_1 B_1 + c_2 B_2 + c_3 B_3 + c_7 B_7 \label{eqnsumB}
  \end{equation}
  we obtain $0 \leqslant B \lesssim \frac{| x |^2}{r} + \frac{| y |^2}{s}$ and
  $\mathd^2 B \geqslant \frac{2}{Q} | \mathd x | | \mathd y |$ in the regions
  where $B \in \mathcal{C}^2$. Indeed,
  \[ (\mathd^2 B_1 \mathd V, \mathd V) \geqslant \frac{4}{Q} | \mathd x | |
     \mathd y | + \frac{4 | x | | y |}{Q r s} | \mathd r | | \mathd s | -
     \frac{4 | y |}{Q s} | \mathd x | | \mathd s | - \frac{4 | x |}{Q r} |
     \mathd y | | \mathd r | \]
  \[ (\mathd^2 B_2 \mathd V, \mathd V) \geqslant \frac{\sqrt{3} | x |}{2 Q r}
     | \mathd y | | \mathd r | - \frac{\sqrt{3} | y |}{2 Q r s} | \mathd r | |
     \mathd s | \]
  \[ (\mathd^2 B_3 \mathd V, \mathd V) \geqslant \frac{\sqrt{3} | y |}{2 Q s}
     | \mathd x | | \mathd s | - \frac{\sqrt{3} | x |}{2 Q r s} | \mathd r | |
     \mathd s | \]
  \[ (\mathd^2 B_7 \mathd V, \mathd V) \geqslant \frac{| x | | y |}{256 Q r s}
     | \mathd r | | \mathd s | \]
  where the last inequality holds in the regions where the function $B_4 \in
  \mathcal{C}^2$. The weighted sum of these inequalities according to
  \ref{eqnsumB} yields the desired inequality on convexity. Now,
  \begin{eqnarray*} 
  && B_1 (V) - B_1 (V_0) - \mathd B_1 (V_0) (V - V_0) \\
  &\gtrsim & \frac{r s}{Q}
     \left| \frac{x}{r} - \frac{x_0}{r_0} \right|  \left| \frac{y}{s} -
     \frac{y_0}{s_0} \right| \geqslant \frac{r s}{Q} \left| \left\langle
     \frac{x}{r} - \frac{x_0}{r_0}, \frac{y}{s} - \frac{y_0}{s_0}
     \right\rangle \right|, 
  \end{eqnarray*}
  \begin{eqnarray*} 
  && B_2 (V) - B_2 (V_0) - \mathd B_2 (V_0) (V - V_0)\\
   & \gtrsim & \frac{s}{Q}
     \frac{| x_0 |}{r_0} | r - r_0 |  \left| \frac{y}{s} - \frac{y_0}{s_0}
     \right| \geqslant \frac{s}{Q} \left| \left\langle \frac{x_0}{r_0},
     \frac{y}{s} - \frac{y_0}{s_0} \right\rangle \right| | r - r_0 |, 
   \end{eqnarray*}
   \begin{eqnarray*} 
   && B_3 (V) - B_3 (V_0) - \mathd B_3 (V_0) (V - V_0) \\
   &\gtrsim &\frac{r}{Q}
     \frac{| y_0 |}{s_0} | s - s_0 | \left| \frac{x}{r} - \frac{x_0}{r_0}
     \right| \geqslant \frac{r}{Q} \left| \left\langle \frac{x}{r} -
     \frac{x_0}{r_0}, \frac{y_0}{s_0} \right\rangle \right| | s - s_0 |, 
   \end{eqnarray*}
   \begin{eqnarray*} 
   && B_7 (V) - B_7 (V_0) - \mathd B_7 (V_0) (V - V_0) \\
   &\gtrsim & \frac{1}{Q} |
     x_0 | | y_0 | | r - r_0 | | s - s_0 | \geqslant \frac{1}{Q} | \langle
     x_0, y_0 \rangle | | r - r_0 | | s - s_0 | . 
    \end{eqnarray*}
  Notice that the last inequalities also remain true when we replace $x$ by
  $\Theta x$ and $x_0$ by $\Theta x_0$ where the rotation $\Theta$ is chosen
  so that $\Theta (x - x_0)$ and $y - y_0$ have the same direction and thus we
  may assume that $\langle x - x_0, y - y_0 \rangle = | x - x_0 | | y - y_0
  |$.
  
  Summing the above inequalities gives
  \begin{eqnarray*}
    &  & Q (B (V) - B (V_0) - \mathd B (V_0) (V - V_0))\\
    & \gtrsim & \left\langle \left( \frac{x}{r} - \frac{x_0}{r_0} \right) r,
    \left( \frac{y}{s} - \frac{y_0}{s_0} \right) s + \frac{y_0}{s_0} (s - s_0)
    \right\rangle\\
    &  & \hspace{1em}+ \left\langle \frac{x_0}{r_0} (r - r_0), \left( \frac{y}{s} -
    \frac{y_0}{s_0} \right) s + \frac{y_0}{s_0} (s - s_0) \right\rangle\\
    & = & \left\langle \left( \frac{x}{r} - \frac{x_0}{r_0} \right) r, y -
    y_0 \right\rangle + \left\langle \frac{x_0}{r_0} (r - r_0), y - y_0
    \right\rangle\\
    & = & \langle x - x_0, y - y_0 \rangle = | x - x_0 | | y - y_0 |
  \end{eqnarray*}
  and we have proved the one-leg convexity. It remains to bound the second
  derivatives in $x$ and $y$. Let $\varepsilon$ be the cut off of the weights
  so that $\varepsilon \leqslant r, s \leqslant \varepsilon^{- 1}$. We
  calculate
  \[ \left( \partial_x^2 \frac{\langle x, x \rangle}{r} \mathd x, \mathd x
     \right) = \frac{2 \langle \mathd x, \mathd x \rangle}{r} \lesssim
     \varepsilon^{- 1} \langle \mathd x, \mathd x \rangle ; \]
  \[ \left( \partial_x^2 \frac{\langle x, x \rangle}{r + M (r, s)} \mathd x,
     \mathd x \right) = \frac{2 \langle \mathd x, \mathd x \rangle}{r + M (r,
     s)} \leqslant \frac{2 \langle \mathd x, \mathd x \rangle}{r} \lesssim
     \varepsilon^{- 1} \langle \mathd x, \mathd x \rangle ; \]
  \begin{eqnarray*}
    &  & \left( \partial_x^2 \frac{\langle x, x \rangle s - 2 | x | | y | K +
    \langle y, y \rangle r}{r s - K^2} \mathd x, \mathd x \right)\\
    & = & \frac{2 \langle \mathd x, \mathd x \rangle s}{r s - K^2} - \frac{2
    | y | K}{r s - K^2} \left( \frac{\langle \mathd x, \mathd x \rangle}{| x
    |} - \frac{\langle x, \mathd x \rangle^2}{| x |^3} \right) \lesssim
    \varepsilon^{- 1} \langle \mathd x, \mathd x \rangle
  \end{eqnarray*}
  where the last implied constant uses the lower bound for $r s - K^2 > r s
  \left( 1 - \frac{1}{Q} \right) \geqslant 1 - \frac{1}{Q}$. We used that
  $\langle x, \mathd x \rangle^2 \leqslant \langle x, x \rangle \langle \mathd
  x, \mathd x \rangle \Rightarrow \frac{\langle \mathd x, \mathd x \rangle}{|
  x |} - \frac{\langle x, \mathd x \rangle^2}{| x |^3} \geqslant 0$.
  
  These imply that for $V \in \mathcal{D}_{Q, \varepsilon}$
  \begin{equation}
    (\partial^2_x B (V) \mathd x, \mathd x) \lesssim \varepsilon^{- 1} |
    \mathd x |^2 . \label{eqnabstau}
  \end{equation}
  This concludes the proof of Lemma \ref{L: existence and properties of the
  Bellman function}.
\end{proof*}

Convexities of the form $\mathd^2 B (V) \geqslant 2 | \mathd x | | \mathd y |$
can be self improved using the following interesting lemma:

\begin{lemma}
  (ellipse lemma, {\tmname{Dragicevic}}--{\tmname{Treil}}--{\tmname{Volberg}}
  {\cite{DraTreVol2008a}}) Let $\mathbbm{H}$ be a Hilbert space with $A, B$
  two positive definite operators on $\mathbbm{H}$. Let $T$ be a self-adjoint
  operator on $\mathbbm{H}$ such that
  \[ (T h, h) \geqslant 2 (A h, h)^{1 / 2} (B h, h)^{1 / 2} \]
  for all $h \in \mathbbm{H}$. Then there exists $\tau > 0$ satisfying
  \[ (T h, h) \geqslant \tau (A h, h) + \tau^{- 1} (B h, h) \]
  for all $h \in \mathbbm{H}$.
\end{lemma}

For our specific Bellman function, we will need a quantitative version:

\begin{lemma}[quantitative ellipse lemma for $B$]
  \label{L: quantitative ellipse Lemma}Let $V \in
  \mathcal{D}^{\varepsilon}_Q$. Assume moreover that $B$ is $\mathcal{C}^2$ at
  $V$. Then there exists $\tau (V) > 0$ such that
  \[ Q \mathd^2_V B (V) \geqslant \tau (V) | \mathd x |^2 + (\tau (V))^{- 1} |
     \mathd y |^2 . \]
  Moreover, we have the bound
  \[ Q^{- 1} \varepsilon \lesssim \tau (V) \lesssim Q \varepsilon^{- 1} .\]
\end{lemma}

\begin{proof*}{Proof of Lemma \ref{L: quantitative ellipse Lemma}}
  {\dueto{quantitative ellipse lemma for $B$}} Let $V \in
  \mathcal{D}^{\varepsilon}_Q$. We have already seen in Lemma \ref{L:
  existence and properties of the Bellman function} that
  \[ \mathd_V^2 B (V) \geqslant \frac{2}{Q} | \mathd x | | \mathd y | . \]
  The ellipse lemma {\cite{DraTreVol2008a}} implies the existence of $\tau
  (V)$ such that for all vectors $\mathd x$ and $\mathd y$ there holds
  \[ Q \mathd^2_V B (V) \geqslant \tau (V) | \mathd x |^2 + (\tau (V))^{- 1}
     | \mathd y |^2 . \]
  We can estimate $\tau (V)$ by testing the Hessian on any $\mathd V$ of the
  form $\mathd V = (\mathd x, 0, 0, 0)$,
  \[ \tau (V) | \mathd x |^2 \leqslant Q (\mathd^2_V B (V) \mathd V, \mathd V)
     = Q (\partial^2_x B (V) \mathd x, \mathd x) \lesssim Q \varepsilon^{- 1}
     | \mathd x |^2 \]
  where the last inequality follows from (\ref{eqnabstau}). Hence $\tau (V)
  \lesssim Q \varepsilon^{- 1}$ as claimed. The same bound holds for $(\tau
  (V)^{- 1})$ by testing against $\mathd V = (0, \mathd y, 0, 0)$. Finally, we
  have proved that for all $V \in \mathcal{D}_Q^{\varepsilon}$,
  \[ Q^{- 1} \varepsilon \lesssim \tau (V) \lesssim Q \varepsilon^{- 1} . \]
\end{proof*}

We now address the lack of smoothness of $B$. All functions aside from $H_4$
that appear are at least in $\mathcal{C}^2$. We apply a standard mollifying
procedure via convolution with $\varphi_{\ell}$ directly on $H_4 (x, y, r, s,
K)$, now only taking {\tmem{real}} variables with $x, y$ positive, $1 < r s <
Q$ and $0 < K < 1$. Here $\varphi$ denotes a standard mollifying kernel in the
five real variables $(x, y, r, s, K) \in \mathbbm{R}^5$ with support in the
corresponding unit ball, whereas $\varphi_{\ell} (\cdot) \assign \ell^{- 5}
\varphi (\cdot / \ell)$ denotes its scaled version with support of size
$\ell$. By slightly changing the constructions, the upper and lower estimate
on the product $r s$ can be modified at the cost of a multiplicative constant
in the final estimate of the Bellman function. Also take into account that the
weights are cut, therefore bounded above and below. Further, we will assume
that the positive variables $x$ and $y$ be bounded below. These considerations
give us enough room to smooth the function $H_4$. It is important that $H_4$
is at least in $\mathcal{C}^1$ and its second order partial derivatives exist
almost everywhere. So we have  $\mathd^2 (H_4 \ast \varphi_{\ell}) = (\mathd^2
H_4) \ast \varphi_{\ell}$. Last, we are observing that as long as the norms of
vectors $| x |$ and $| y |$ are bounded away from 0, our function $H_4 \ast
\varphi_{\ell}$, mollified in $\mathbbm{R}^5$ remains smooth when taking
vector variables (observe that the final Bellman function only depends upon $|
x |$ and $| y |$). It is important that the smoothing happens before the
function is composed with $K$, we therefore preserve fine convexity
properties, in particular also the much needed one-leg convexity. Size
estimates change slightly, but are recovered when the mollifying parameter
goes to $0$. These details are either standard and have appeared in numerous
articles on Bellman functions or an easy consequence of reading the
construction of the Bellman function above.

\begin{lemma}[regularised Bellman function and its properties]
  \label{L: regularised Bellman function}Let $\varepsilon > 0$ given. Let $0 <
  \ell \leqslant \varepsilon / 2$. There exists a function $B_{\ell} (x, y, r,
  s)$ defined with domain
  \[ \mathcal{D}_Q^{\varepsilon, \ell} \assign \left\{ V \in
     \mathcal{D}_Q^{\varepsilon} ; \hspace{1em} | x | \geqslant \ell, | y |
     \geqslant \ell \right\} \subset \mathcal{D}_Q^{\varepsilon} \]
  such that for all $V_0, V \in \mathcal{D}_Q^{\varepsilon, \ell}$, we have
  \[ \hspace{1em} B_{\ell} \lesssim (1 + \ell) \left( \frac{| x |^2}{r} +
     \frac{| y |^2}{s} \right), \]
  \begin{equation}
    \mathd_V^2 B_{\ell} (V) \geqslant \frac{2}{Q}  | \mathd x | | \mathd y |
    \nocomma, \label{eq: concavity continuous}
  \end{equation}
  \begin{equation}
    B_{\ell} (V) - B_{\ell} (V_0) - \mathd_V B_{\ell} (V_0) (V - V_0)
    \geqslant \frac{1}{Q} | \Delta x |  | \Delta y | = \frac{1}{Q} | x - x_0 |
    | y - y_0 | \label{eq: concavity discrete}
  \end{equation}
  and moreover the quantitative ellipse lemma now holds in the form
  \begin{eqnarray*}
  && Q \mathd^2_V B_{\ell} (V)\\
  & \geqslant &\tau_{\ell} (V) | \mathd x |^2 +
     (\tau_{\ell} (V))^{- 1} | \mathd y |^2 , 
   \end{eqnarray*}
  where $\tau_{\ell} \assign \tau_{\ell} (V)$ is a continuous function of its
  arguments, and where
  \[ Q^{- 1} \varepsilon \lesssim \tau_{\ell} (V) \lesssim Q \varepsilon^{- 1}
     . \]
\end{lemma}

\

\section{Dissipation estimates}

\

Let $V \assign (X, Z, u, w)$ a c{\`a}dl{\`a}g adapted martingale with values
in $\mathcal{D}_Q^{\varepsilon}$. In order to bound away from the
$\mathbbm{H}$--valued martingale $X \assign (X^1, X^2, \ldots)$, it is
classical to introduce the $\mathbbm{R} \times \mathbbm{H}$--valued
martingales $X^a \assign (a, X^1, X^2, \ldots)$ where $a > 0$. It follows that
$\| X^a \|^2 = \| X \|^2 + a^2$ and $\| X^a \| \geqslant a$, and the same
construction holds for $Z$. We note $V^a \assign (X^a, Z^a, u, w)$ Given $\ell
> 0$ a smoothing parameter, take $a \geqslant \ell$ then it follows that
\[ V \in \mathcal{D}_Q^{\varepsilon} \Rightarrow V^a \in
   \mathcal{D}_Q^{\varepsilon, \ell} .\]
The main result of this section is the following dissipation estimate:

\begin{proposition}[dissipation estimate]
  \label{P: Dissipation estimates}Let $\varepsilon > 0$, $\ell > 0$ as defined
  above. Let $V$ a c{\`a}dl{\`a}g adapted martingale with $V \in
  \mathcal{D}^{\varepsilon}_Q$. Let \ $F_t \assign \mathbbm{E} (| X_{\infty}
  |^2 w^{\varepsilon}_{\infty}  |  \mathcal{F}_t)$ and $G_t \assign
  \mathbbm{E} (| Z_{\infty} |^2 u^{\varepsilon}_{\infty}  |  \mathcal{F}_t)$.
  Let finally $a \geqslant \ell$. We have
  \begin{eqnarray*}
    && Q (1 + \ell) (\mathbbm{E}F_t +\mathbbm{E}G_t + 2 a^2 \varepsilon^{- 1}) \\
    &
    \gtrsim & \frac{1}{2} \mathbbm{E} \int_0^t \tau_{\ell} (V_{s -}) \mathd
    [X, X]^c_s + (\tau_{\ell} (V_{s -}))^{- 1} \mathd [Z, Z]^c_s\\
    &  & \hspace{2em} +\mathbbm{E} \sum_{0 < s \leqslant t} | \Delta X_s |  |
    \Delta Z_s | .
  \end{eqnarray*}
\end{proposition}

We need the preliminary lemma

\begin{lemma}[comparison of quadratic forms in stochastic integrals]
  \label{L: comparison of quadratic forms}Let $\mathcal{Q}$ denote the set of
  quadratic forms from $\mathbbm{R}^m \times \mathbbm{R}^m \rightarrow
  \mathbbm{R}$. Let $A \assign (A_{\alpha \beta})_{1 \leqslant \alpha, \beta
  \leqslant m}$ and $B \assign (B_{\alpha \beta})_{1 \leqslant \alpha, \beta
  \leqslant m}$ two $\mathcal{Q}$--valued c{\`a}dl{\`a}g processes. Assume for
  all $t \geqslant 0$ and a.s. that $A (t) \geqslant B (t)$ (resp. $A (t)
  \geqslant | B (t) |$), in the sense $\forall \mathd V \in \mathbbm{R}^m,$
  \[ \hspace{1em} (A \mathd V, \mathd V)
     \geqslant (B \mathd V, \mathd V) \hspace{1em} (\tmop{resp} . (A \mathd V,
     \mathd V) \geqslant | (B \mathd V, \mathd V) |) . \]
  Abbreviating $A_{s -} : \mathd [V, V]_s = \sum_{\alpha, \beta} (A_{\alpha
  \beta})_{s -} \mathd [V_{\alpha}, V_{\beta}]_s$ then for all $t \geqslant
  0$,
  \[ \mathbbm{E} \int_0^t A_{s -} : \mathd [V, V]_s \geqslant \mathbbm{E}
     \int_0^t B_{s -} : \mathd [V, V]_s, \]
  \[ \left( \tmop{resp} .\mathbbm{E} \int_0^t A_{s -} : \mathd [V, V]_s
     \geqslant \mathbbm{E} \int_0^t | B_{s -} : \mathd [V, V]_s | \right) . \]
\end{lemma}

\

\begin{proof*}{Proof of Lemma \ref{L: comparison of quadratic forms}}
  {\dueto{comparison of quadratic forms in stochastic integrals}}With the
  hypotheses above, let us consider the case $A (t) \geqslant B (t)$, the case
  $A (t) \geqslant | B (t) |$ being treated in the same manner. Given $t
  \geqslant 0$, assume that
  \[ \int_0^t A_{s -} : \mathd [V, V]_s = \sum_{\alpha, \beta} \int_0^t
     (A_{\alpha \beta})_{s -} \mathd [V_{\alpha}, V_{\beta}]_s < \infty \]
  otherwise the claim is proved. Given the process $V$, let $\sigma_n \assign
  (0 \leqslant T_0^n \leqslant T_1^n \leqslant \ldots \leqslant T_i^n
  \leqslant \ldots \leqslant T_{k_n}^n \leqslant t)$ denote a random partition
  of stopping times tending to the identity as $n$ tends to infinity. Given
  $\alpha$ and $\beta$, we have that $A_{\alpha \beta}$ is a
  $\mathbbm{R}$--valued c{\`a}dl{\`a}g process. It follows (see e.g. Protter
  {\cite{Pro2005a}}) that the stochastic integral
  \begin{equation}
    \int_0^t A_{\alpha \beta} (s -) \mathd [V_{\alpha}, V_{\beta}]_s
    \label{eq: stochastic integral Aalphabeta}
  \end{equation}
  is the limit in \tmtextit{ucp} (uniform convergence in probability) as $n$
  tends to infinity of sums
  \[ S^A_{\alpha \beta} \assign \sum_{i = 0}^{k_n - 1} A_{\alpha \beta}
     (T^n_i) (V_{\alpha}^{T^n_{i + 1}} - V_{\alpha}^{T^n_i})
     (V_{\beta}^{T^n_{i + 1}} - V_{\beta}^{T^n_i}) \]
  involving the stopping times defined above. Since $A \geqslant B$, summing
  w.r.t. $\alpha, \beta$ yields, for any $s \in [0, t]$,
  \begin{eqnarray*}
    \left( \sum_{\alpha, \beta} S^A_{\alpha \beta} \right) (s) & \assign &
    \sum_{\alpha, \beta} \sum_{i = 0}^{k_n - 1} A_{\alpha \beta} (T^n_i)
    (V_{\alpha, s}^{T^n_{i + 1}} - V_{\alpha, s}^{T^n_i}) (V_{\beta,
    s}^{T^n_{i + 1}} - V_{\beta, s}^{T^n_i})\\
    & = & \sum_{i = 0}^{k_n - 1} \sum_{\alpha, \beta} A_{\alpha \beta}
    (T^n_i) (V_{\alpha, s}^{T^n_{i + 1}} - V_{\alpha, s}^{T^n_i}) (V_{\beta,
    s}^{T^n_{i + 1}} - V_{\beta, s}^{T^n_i})\\
    & \geqslant & \sum_{i = 0}^{k_n - 1} \sum_{\alpha, \beta} B_{\alpha
    \beta} (T^n_i) (V_{\alpha, s}^{T^n_{i + 1}} - V_{\alpha, s}^{T^n_i})
    (V_{\beta, s}^{T^n_{i + 1}} - V_{\beta, s}^{T^n_i})\\
    & \geqslant & \left( \sum_{\alpha, \beta} S^B_{\alpha \beta} \right) (s)
  \end{eqnarray*}
  with an obvious definition for $S^B_{\alpha \beta}$. Passing to the limit in
  the sums $\sum_{\alpha, \beta}$ gives the result.
\end{proof*}

\

\begin{proof*}{Proof of Proposition \ref{P: Dissipation estimates}}
  {\dueto{dissipation estimates}}
  
  \paragraph{Step 1}We first pas to a finite dimensional case. Let $V$ a
  c{\`a}dl{\`a}g adapted martingale with $V \in \mathcal{D}^{\varepsilon}_Q$.
  Then $V^a \in \mathcal{D}_Q^{\varepsilon, \ell}$. We note $X^{a, m}$ the
  projection of $X^a \in \mathbbm{R} \times \mathbbm{H}$ onto $\mathbbm{R}
  \times \mathbbm{R}^m$, and introduce accordingly $Z^{a, m}$ and $V^{a, m}$.
  Notice that $[X^a, X^a] = a^2 + [X, X]$ and similarly $[X^{a, m}, X^{a, m}]
  = a^2 + [X^m, X^m]$. Since $V^{a, m} \in \mathcal{D}_Q^{\varepsilon, \ell}$
  where $B_{\ell}$ is $\mathcal{C}^2$ and we can apply It{\^o}'s formula and
  obtain, for all $t > 0$, almost sure paths,
  \begin{eqnarray*}
    && B_{\ell} (V_t^{a, m}) \\
    & = & B_{\ell} (V_0^{a, m}) + \int_{0 +}^t \mathd_V
    B (V_{s -}^{a, m}) \mathd V^m_s + \frac{1}{2} \int_{0 +}^t \mathd^2_V
    B_{\ell} (V^{a, m}_{s -}) : \mathd [V^m, V^m]^c_s\\
    &  & \hspace{2em} + \sum_{0 < s \leqslant t} \{ B_{\ell} (V^{a, m}_s) -
    B_{\ell} (V^{a, m}_{s -}) - \mathd_V B_{\ell} (V^{a, m}_{s -}) \Delta
    V^m_s \} .
  \end{eqnarray*}
  Thanks to Lemma \ref{L: quantitative ellipse Lemma} and Lemma \ref{L:
  comparison of quadratic forms}, the concavity properties (\ref{eq: concavity
  continuous}) of $B_{\ell}$ imply for the continuous part
  \begin{eqnarray*}
  && \frac{1}{2} \int_{0 +}^t \mathd^2_V B_{\ell} (V^{a, m}_{s -}) : \mathd
     [V^m, V^m]^c_s \\
     & \geqslant &
     \frac{1}{2 Q} \int_{0 +}^t \tau_{\ell} (V^{a,
     m}_{s -}) \mathd [X^m, X^m]^c_s + (\tau_{\ell} (V^{a, m}_{s -}))^{- 1}
     \mathd [Z^m, Z^m]^c_s . 
     \end{eqnarray*}
  Also, the concavity properties (\ref{eq: concavity discrete}) of $B_{\ell}$
  for the jump part
  \[ B_{\ell} (V^{a, m}_s) - B_{\ell} (V^{a, m}_{s -}) - \mathd_V B_{\ell}
     (V^{a, m}_{s -}) \Delta V^m_s \geqslant \frac{1}{Q} | \Delta X^m_s |  |
     \Delta Z^m_s | . \]
  Plugging the continuous and jump dissipation estimates into It{\^o}'s
  formula yields for all times, almost sure paths,
  \begin{eqnarray*}
    B_{\ell} (V^{a, m}_t) & \geqslant & B_{\ell} (V^{a, m}_0) + \int_{0 +}^t
    \mathd_V B_{\ell} (V^{a, m}_{s -}) \mathd V^m_s\\
    &  & + \frac{1}{2 Q} \int_0^t \tau_{\ell} (V^{a, m}_{s -})
    \mathd [X^m, X^m]^c_s + (\tau_{\ell} (V^m_{s -}))^{- 1} \mathd [Z^m,
    Z^m]^c_s\\
    &  & + \frac{1}{Q} \sum_{0 < s \leqslant t} | \Delta X^m_s |
    | \Delta Z^m_s | .
  \end{eqnarray*}

  \paragraph{Step 2}For technical reasons in the proof, we work with
  bounded martingales that we obtain through a usual stopping procedure.
  Recall that $V$ is a c{\`a}dl{\`a}g adapted martingale with $V \in
  \mathcal{D}^{\varepsilon}_Q$ and $V^a \in \mathcal{D}^{\varepsilon,
  \ell}_Q$. For all $M \in \mathbbm{N}$, define the stopping time $T_M$ as
  $T_M \assign \inf \{ t > 0 ; | V^a |^2_t + [V^a, V^a]_t > M^2 \}$, so that
  $T_M$ is a stopping time that tends to infinity as $M$ goes to infinity. It
  follows that $V^{a, T_M}$ is a local martingale, and that $V^{a, T_M -}$ and
  $[V^a, V^a]^{T_M -}$ are bounded semimartingales. Let $m \in
  \mathbbm{N}^{\star}$ and $V^{a, m}$ the projection of $V^a$ onto
  $\mathbbm{R}^m \subset \mathbbm{H}$. For each $M$, there exists a sequence
  $\{ T_{M, k} \}_{k \geqslant 1}$ of stopping times such that $T_{M, k}
  \hspace{1em} \uparrow \hspace{1em} T_M$ as $k \uparrow \infty$, and such
  that $(V^{a, m})^{T_{M, k}}$ is a martingale. Since $| V^{a, m} | \leqslant
  | V^a |$, it follows that $(V^{a, m})^{T_{M, k} -}$ is a bounded
  semimartingale, to which we can apply the dissipation estimate of Step 1
  above and obtain
  \begin{eqnarray*}
    && B_{\ell} (V^{a, m}_{t \wedge T_{M, k} -}) \\
    & \geqslant & B_{\ell} (V^{a,
    m}_0) + \int_{0 +}^{t \wedge T_{M, k} -} \mathd_V B_{\ell} (V^{a, m}_{s
    -}) \mathd V^m_s\\
    &  & \hspace{1em} + \frac{1}{2 Q} \int_0^{t \wedge T_{M, k} -}
    \tau_{\ell} (V^{a, m}_{s -}) \mathd [X^m, X^m]^c_s + (\tau_{\ell} (V^{a,
    m}_{s -}))^{- 1} \mathd [Z^m, Z^m]^c_s\\
    &  & \hspace{1em} + \frac{1}{Q} \sum_{0 < s < t \wedge T_{M, k}} | \Delta
    X^m_s |  | \Delta Z^m_s |\\
    & = & B_{\ell} (V^{a, m}_0) + \int_{0 +}^{t \wedge T_{M, k}} \mathd_V
    B_{\ell} (V^{a, m}_{s -}) \mathd V^m_s\\
    &  & \hspace{1em} + \frac{1}{2 Q} \int_0^{t \wedge T_{M, k} -}
    \tau_{\ell} (V^{a, m}_{s -}) \mathd [X^m, X^m]^c_s + (\tau_{\ell} (V^{a,
    m}_{s -}))^{- 1} \mathd [Z^m, Z^m]^c_s\\
    &  & \hspace{1em} + \frac{1}{Q} \sum_{0 < s < t \wedge T_{M, k}} | \Delta
    X^m_s |  | \Delta Z^m_s | - \mathd_V B_{\ell} (V^{a, m}_{t \wedge T_{M, k} -})
    \Delta V^m_{t \wedge T_{M, k}} .
  \end{eqnarray*}
  Taking expectation and then letting $k \rightarrow \infty$, the dominated
  convergence theorem yields
     \begin{eqnarray}
     &&  \mathbbm{E}B_{\ell} (V^{a, m}_{t \wedge T_M -})  \label{i:withmstill} \\
    & \geqslant &
    \mathbbm{E}B_{\ell} (V^{a, m}_0)\nonumber \\
    &&   \hspace{1em}+ \frac{1}{2 Q} \mathbbm{E} \int_0^{t
    \wedge T_M} \tau_{\ell} (V^{a, m}_{s -}) \mathd [X^m, X^m]^c_s +
    (\tau_{\ell} (V^{a, m}_{s -}))^{- 1} \mathd [Z^m, Z^m]^c_s \nonumber\\
    &  & \hspace{1em} + \frac{1}{Q} \mathbbm{E} \sum_{0 < s < t \wedge T_M} |
    \Delta X^m_s |  | \Delta Z^m_s | \nonumber -\mathbbm{E} \{ \mathd_V B_{\ell} (V^{a, m}_{t \wedge
    T_M -}) \Delta V^m_{t \wedge T_M} \nonumber \} .  
  \end{eqnarray}
  Observe that we used size properties of $B_{\ell}$, the definition of the
  stopping time $T_{M, k}$, the estimate of the $\tau_{\ell}$ provided by
  Lemma \ref{L: regularised Bellman function} and the size control of the
  weights.
  
  \
  
  \paragraph{Step 3}Now, we wish to return to the infinite dimensional case.
  First recall that $0 \leqslant B_{\ell} (V) \lesssim (1 + \ell) (X^2 / u +
  Y^2 / w)$. 
  
  Let $F_t \assign \mathbbm{E} (| X_{\infty} |^2 w_{\infty} |
  \mathcal{F}_t)$, $G_t \assign \mathbbm{E} (| Z_{\infty} |^2 u_{\infty} |
  \mathcal{F}_t)$ as well as $F^a_t \assign \mathbbm{E} (| X^a_{\infty} |^2 w_{\infty} |
  \mathcal{F}_t)$, $G^a_t \assign \mathbbm{E} (| Z^a_{\infty} |^2
  u_{\infty}  | \mathcal{F}_t)$. Notice that $F^a_t =\mathbbm{E} ((|
  X_{\infty} |^2 + a^2) w_{\infty} | \mathcal{F}_t) \leqslant F_t +\mathbbm{E}
  (a^2 w_{\infty}^{\varepsilon} | \mathcal{F}_t) \leqslant F_t + a^2
  \varepsilon^{- 1}$. It follows, thanks to Jensen inequality, that
  \[ B_{\ell} (V_t^a) \leqslant C_0  (1 + \ell) (F^a_t + G^a_t) \lesssim (1 +
     \ell) (F_t + G_t + 2 a^2 \varepsilon^{- 1}) \]
  A similar inequality holds for $V^{a, m}$ and in particular
  \begin{eqnarray*}
    B_{\ell} (V^{a, m}_{t \wedge T_M -}) & \lesssim &  (1 + \ell) (F_{t \wedge
    T_M -} + G_{t \wedge T_M -} + 2 a^2 \varepsilon^{- 1}) .
  \end{eqnarray*}
  Hence, the dominated convergence theorem implies that $\mathbbm{E}B_{\ell}
  (V^{a, m}_{t \wedge T_M -})$ converges when $m$ goes to infinity towards
  $\mathbbm{E}B_{\ell} (V^a_{t \wedge T_M -})$.
  
  \
  
  Let us consider the first term in the last integral of step 2, the second
  term integral in inequality (\ref{i:withmstill}). We write
  \begin{eqnarray*}
    &&\mathbbm{E} \int_0^{t \wedge T_M} \tau_{\ell} (V^{a, m}_{s -}) \mathd
    [X^m, X^m]^c_s \\
    & = & \mathbbm{E} \int_0^{t \wedge T_M} \tau_{\ell} (V^a_{s
    -}) \mathd [X, X]^c_s\\
    &  & \hspace{2em} +\mathbbm{E} \int_0^{t \wedge T_M} (\tau_{\ell} (V^{a,
    m}_{s -}) - \tau_{\ell} (V^a_{s -})) \mathd [X, X]^c_s\\
    &  & \hspace{2em} +\mathbbm{E} \int_0^{t \wedge T_M} \tau_{\ell} (V^{a,
    m}_{s -}) \mathd ([X^m, X^m]^c - [X, X]^c)_s .
  \end{eqnarray*}
  The uniform boundedness and continuity of $\tau_{\ell}$, the square
  integrability of $X$ and the Dominated convergence theorem imply that the
  second term of the right--hand--side converges to zero. The last term can be
  bounded above using the estimates for $\tau_{\ell}$.
  \begin{eqnarray*}
  && \left| \mathbbm{E} \int_0^{t \wedge T_M} \tau_{\ell} (V^{a, m}_{s -})
    \mathd ([X^m, X^m]^c - [X, X]^c)_s \right| \\
    & \lesssim &
    \frac{Q}{\varepsilon} \mathbbm{E} \int_0^{t \wedge T_M}  | \mathd ([X^m,
    X^m]^c - [X, X]^c)_s |\\
    & \lesssim & \frac{Q}{\varepsilon} \mathbbm{E} \int_0^{t \wedge T_M}
    \mathd ([X, X]^c - [X^m, X^m]^c)_s\\
    & \lesssim & \frac{Q}{\varepsilon}  (\mathbbm{E} [X, X]_{t \wedge T_M}^c
    -\mathbbm{E} [X^m, X^m]_{t \wedge T_M}^c)
  \end{eqnarray*}
  where we used that for $m$ fixed, $[X, X]^c - [X^m, X^m]^c$ is a nonnegative
  nondecreasing process. The last expression in the last line tends to zero
  when $m \rightarrow \infty$ by the monotone convergence theorem. We prove in
  a similar manner the convergence
  \[ \mathbbm{E} \sum_{0 < s < t \wedge T_M} | \Delta X^m_s |  | \Delta Z^m_s
     | \xrightarrow[m \rightarrow \infty]{} \mathbbm{E} \sum_{0 < s < t \wedge
     T_M} | \Delta X_s |  | \Delta Z_s | . \]

  Finally, since $| V^{a, m}_{t \wedge T_M -} | \leqslant | V^a_{t \wedge T_M
  -} |$ for all $m$, $\mathd_V B_{\ell}$ is continuous and bounded on
  compacts, $| \Delta V^m_{t \wedge T_M} |^2 \leqslant | \Delta V_{t \wedge
  T_M} |^2 \leqslant [V, V]_t$ and $\mathbbm{E} [V, V]_t =\mathbbm{E} | V_t
  |^2 < \infty$, the dominated convergence theorem ensures that
  \[ -\mathbbm{E} \{ \mathd_V B_{\ell} (V^{a, m}_{t \wedge T_M -}) \Delta
     V^m_{t \wedge T_M} \} \rightarrow -\mathbbm{E} \{ \mathd_V B_{\ell}
     (V^a_{t \wedge T_M -}) \Delta V_{t \wedge T_M} \} . \]
  Collecting all terms,
  \begin{eqnarray*}
    &&\mathbbm{E}B_{\ell} (V^a_{t \wedge T_M -})\\
     & \geqslant & \frac{1}{2 Q}
    \mathbbm{E} \int_0^{t \wedge T_M} \tau_{\ell} (V^a_{s -}) \mathd [X,
    X]^c_s + (\tau_{\ell} (V^a_{s -}))^{- 1} \mathd [Z, Z]^c_s\\
    &  & \hspace{2em} + \frac{1}{Q} \mathbbm{E} \sum_{0 < s < t \wedge T_M} |
    \Delta X_s |  | \Delta Z_s |\\
    &  & \hspace{2em} -\mathbbm{E} \{ \mathd_V B_{\ell} (V^a_{t \wedge T_M
    -}) \Delta V_{t \wedge T_M} \} .
  \end{eqnarray*}

  \paragraph{Step 4}Now we add the contribution of the possible jumps occuring
  at $T_M$. We have seen in Step 1 the dissipation estimate along one jump
  \begin{eqnarray*}
   &&B_{\ell} (V^a_{t \wedge T_M}) - B_{\ell} (V^a_{t \wedge T_M -}) -
     \mathd_V B_{\ell} (V^a_{t \wedge T_M -}) \Delta V_{t \wedge T_M} \\
     &\geqslant & \frac{1}{Q}  | \Delta X_{t \wedge T_M} |  | \Delta Z_{t \wedge
     T_M} | . \end{eqnarray*}
  Taking expectation and adding the contribution of Step 3 yields
  \begin{eqnarray}
    \mathbbm{E}B_{\ell} (V^a_{t \wedge T_M}) & \geqslant & \frac{1}{2 Q}
    \mathbbm{E} \int_0^{t \wedge T_M} \tau_{\ell} (V^a_{s -}) \mathd [X,
    X]^c_s + (\tau_{\ell} (V^a_{s -}))^{- 1} \mathd [Z, Z]^c_s \nonumber\\
    &  & \hspace{2em} + \frac{1}{Q} \mathbbm{E} \sum_{0 < s \leqslant t
    \wedge T_M} | \Delta X_s |  | \Delta Z_s | .  \label{i:withMstill}
  \end{eqnarray}

  \paragraph{Step 5}We will pass to the limit $M \rightarrow \infty$. Recall
  again that $0 \leqslant B_{\ell} (V) \lesssim (1 + \ell) (X^2 / u + Y^2 /
  w)$. Using Doob's inequality for square integrable martingales, we have for
  all $M$
  \[ \mathbbm{E}B_{\ell} (V^a_{t \wedge T_M}) \lesssim \varepsilon^{- 1} (1 +
     \ell) (\mathbbm{E}X^2 +\mathbbm{E}Y^2) < \infty . \]
  So, by the dominated convergence theorem, $\mathbbm{E}B_{\ell} (V^a_{t
  \wedge T_M}) \rightarrow \mathbbm{E}B_{\ell} (V^a_t)$ as $M \rightarrow
  \infty$. The monotone convergence theorem for the integral in the
  right--hand--side of the inequality (\ref{i:withMstill}) therefore yields in
  the limit $M \rightarrow \infty$
  \begin{eqnarray*}
    &&(1 + \ell) (\mathbbm{E}F_t +\mathbbm{E}G_t + 2 a^2 \varepsilon^{- 1}) \\
    & \geqslant & \mathbbm{E}B_{\ell} (V^a_t)\\
    & \gtrsim & \frac{1}{2 Q} \mathbbm{E} \int_0^t \tau_{\ell} (V^a_{s -})
    \mathd [X, X]^c_s + (\tau_{\ell} (V^a_{s -}))^{- 1} \mathd [Z, Z]^c_s\\
    &  & \hspace{2em} + \frac{1}{Q} \mathbbm{E} \sum_{0 < s \leqslant t} |
    \Delta X_s |  | \Delta Z_s | .
  \end{eqnarray*}
  This concludes the proof of Proposition \ref{P: Dissipation estimates}.
  
  \ 
\end{proof*}

\section{Truncation of the weights}

\

Due to several technicalities in the proof, we have used weights bounded from
above and away from 0. In order to pass to the general case, we cut a possibly
unbounded weight above and below and show that this operation does not
increase the charateristic of the weight. This is convenient and has been used
in several places, here we extend {\cite{RezVasVol2010a}} to the martingale
setting. Their proof is particularly nice, since it does not increase the
characteristic at all, not even by a constant. We need the following
preliminary lemmas.

\begin{lemma}[truncation from above]
  \label{L: truncation above}For $a > 0$ let $M = \{ w \leqslant a \}$ and $H
  = \{ w > a \}$. Now take $w_{\bar{a}} = w \chi_M + a \chi_H$. Then
  $Q^{\mathcal{F}}_2 [w_{\bar{a}}] \leqslant Q^{\mathcal{F}}_2 [w]$.
\end{lemma}

\begin{proof*}{Proof of Lemma \ref{L: truncation above}}
  Let $\tau$ be a stopping time and let us decompose
  \begin{eqnarray*}
   \mathbbm{E} (w | \mathcal{F}_{\tau}) & = & \mathbbm{E} (w \chi_M |
     \mathcal{F}_{\tau}) +\mathbbm{E} (w \chi_H | \mathcal{F}_{\tau})\\
     &=&\mathbbm{E} (\chi_M | \mathcal{F}_{\tau}) \mathbbm{E}_M (w |
     \mathcal{F}_{\tau}) +\mathbbm{E} (\chi_H | \mathcal{F}_{\tau})
     \mathbbm{E}_H (w | \mathcal{F}_{\tau}) \end{eqnarray*}
  where for example $\mathbbm{E}_M (w | \mathcal{F}_{\tau})$ means expectation
  is taken with respect to the measure $\chi_M \mathd \mathbbm{P}$. Write as
  usual $\mathbbm{E} (\chi_M | \mathcal{F}_{\tau}) = (\chi_M)_{\tau}$.
  \begin{eqnarray*}
    &  & \mathbbm{E} (w | \mathcal{F}_{\tau}) \mathbbm{E} (w^{- 1} |
    \mathcal{F}_{\tau}) -\mathbbm{E} (w_{\bar{a}} | \mathcal{F}_{\tau})
    \mathbbm{E} (w_{\bar{a}}^{- 1} | \mathcal{F}_{\tau})\\
    & = & ((\chi_M)_{\tau} \mathbbm{E}_L (w | \mathcal{F}_{\tau}) +
    (\chi_H)_{\tau} \mathbbm{E}_H (w | \mathcal{F}_{\tau})) ((\chi_M)_{\tau}
    \mathbbm{E}_L (w^{- 1} | \mathcal{F}_{\tau}) \\
    &&\hspace{2em} + (\chi_H)_{\tau}
    \mathbbm{E}_H (w^{- 1} | \mathcal{F}_{\tau})) - ((\chi_M)_{\tau} \mathbbm{E}_L (w | \mathcal{F}_{\tau}) \\
    && \hspace{2em}+
    (\chi_H)_{\tau} a) ((\chi_M)_{\tau} \mathbbm{E}_L (w^{- 1} |
    \mathcal{F}_{\tau}) + (\chi_H)_{\tau} a^{- 1})\\
    & = & (\chi_M)_{\tau} (\chi_H)_{\tau} \left( \mathbbm{E}_M (w |
    \mathcal{F}_{\tau}) \mathbbm{E}_H (w^{- 1} | \mathcal{F}_{\tau})
    +\mathbbm{E}_M (w^{- 1} | \mathcal{F}_{\tau}) \mathbbm{E}_H (w |
    \mathcal{F}_{\tau}) \right. \\
    &&\hspace{5em}\left.-\mathbbm{E}_M (w | \mathcal{F}_{\tau}) a^{- 1}
    -\mathbbm{E}_M (w^{- 1} | \mathcal{F}_{\tau}) a \right)\\
    &  & \hspace{2em}+ (\chi_H)^2_{\tau} (\mathbbm{E}_H (w | \mathcal{F}_{\tau})
    \mathbbm{E}_H (w^{- 1} | \mathcal{F}_{\tau}) - 1) .
  \end{eqnarray*}
  The last term is positive thanks to Jensen inequality. Let us observe that
  also
  \begin{eqnarray*}
    &  & \mathbbm{E}_M (w | \mathcal{F}_{\tau}) \mathbbm{E}_H (w^{- 1} |
    \mathcal{F}_{\tau}) +\mathbbm{E}_M (w^{- 1} | \mathcal{F}_{\tau})
    \mathbbm{E}_H (w | \mathcal{F}_{\tau}) \\
    &&\hspace{2em}-\mathbbm{E}_M (w |
    \mathcal{F}_{\tau}) a^{- 1} -\mathbbm{E}_M (w^{- 1} | \mathcal{F}_{\tau})
    a\\
    & = & \mathbbm{E}_M (w | \mathcal{F}_{\tau}) \mathbbm{E}_H (w^{- 1} -
    a^{- 1} | \mathcal{F}_{\tau}) +\mathbbm{E}_M (w^{- 1} |
    \mathcal{F}_{\tau}) \mathbbm{E}_H (w - a | \mathcal{F}_{\tau})\\
    & = & \mathbbm{E}_H \left( \frac{w - a}{w a} (w a\mathbbm{E}_M (w^{- 1} |
    \mathcal{F}_{\tau}) -\mathbbm{E}_M (w | \mathcal{F}_{\tau})) \middle|
    \mathcal{F}_{\tau} \right)\\
    & \geqslant & 0.
  \end{eqnarray*}
  Here the last inequality uses $\mathbbm{E}_M (w^{- 1} | \mathcal{F}_{\tau})
  \geqslant a^{- 1}$ and $\mathbbm{E}_M (w | \mathcal{F}_{\tau}) \leqslant a$
  also $w - a \geqslant 0$ on $H$. This proves the Lemma.
\end{proof*}

\begin{lemma}[two-sided truncation]
  \label{L: truncation above and below}For $a > 0$ let $M = \{ a^{- 1}
  \leqslant w \leqslant a \}$ and $L = \{ w < a^{- 1} \}$ and $H = \{ w > a
  \}$ then with $w_a = a^{- 1} \chi_L + w \chi_M + a \chi_H$ we have
  $Q^{\mathcal{F}}_2 [w_a] \leqslant Q^{\mathcal{F}}_2 [w] .$
\end{lemma}

\begin{proof*}{Proof of Lemma \ref{L: truncation above and below}}
  Let $w_{\bar{a}}$ be the weight obtained in the previous lemma. Apply now
  the previous lemma to $w_{\bar{a}}^{- 1}$, truncating above by the same $a$.
\end{proof*}

\

\section{Proof of the main results}

\

\begin{proof*}{Proof of Proposition \ref{P: bilinear estimate}}
  {\dueto{bilinear estimate}}Let $\lambda > 0$. Let $Y$ differentially
  subordinate to $X$, then $\lambda Y$ is differentially subordinate to
  $\lambda X$. Let $w$ a weight in the $\tmmathbf{A}_2$ class. Let
  $w^{\varepsilon}$ the $\varepsilon$--truncation of $w$. Using Proposition
  \ref{P: Dissipation estimates} with $V^{\varepsilon, \lambda} \assign
  (\lambda X, \lambda^{- 1} Z, u^{\varepsilon}, w^{\varepsilon})$ and $Q =
  Q^{\mathcal{F}}_2 [w]$. Notice that since $Q^{\mathcal{F}}_2
  [w^{\varepsilon}] \leqslant Q^{\mathcal{F}}_2 [w]$, $V^{\varepsilon,
  \lambda} \in \mathcal{D}^{\varepsilon, \ell}_Q$ using the differential
  subordination of $\lambda Y$ w.r.t. $\lambda X$, we have for all $t > 0$,
   \begin{eqnarray*}
       &&Q^{\mathcal{F}}_2 [w] (1 + \ell) (\mathbbm{E} \lambda^2 F_t
       +\mathbbm{E} \lambda^{- 2} G_t + 2 a^2 \varepsilon^{- 1})\\
       &\gtrsim & \frac{1}{2} \mathbbm{E} \int_0^t \tau_{\ell}
       (V^a_{s -}) \mathd [\lambda X, \lambda X]^c_s + (\tau_{\ell} (V^a_{s
       -}))^{- 1} \mathd [\lambda^{- 1} Z, \lambda^{- 1} Z]^c_s\\
       &&\hspace{2em}+\mathbbm{E} \sum_{0 < s \leqslant t} | \lambda \Delta X_s
       |  | \lambda^{- 1} \Delta Z_s |\\
       &\gtrsim & \frac{1}{2} \mathbbm{E} \int_0^t \tau_{\ell}
       (V^a_{s -}) \mathd [\lambda Y, \lambda Y]^c_s + (\tau_{\ell} (V_{s
       -}))^{- 1} \mathd [\lambda^{- 1} Z, \lambda^{- 1} Z]^c_s\\
        &&\hspace{2em}+\mathbbm{E} \sum_{0 < s \leqslant t} | \Delta Y_s |  |
       \Delta Z_s | .
     \end{eqnarray*} 
  Since for any $0 < \kappa < \infty$ and any $x \in \mathbbm{H}$, $y \in
  \mathbbm{H}$, we have $\kappa x^2 + \kappa^{- 1} y^2 \geqslant 2 | \langle
  x, y \rangle |$, it follows easily
  \begin{eqnarray*}
       &&\frac{1}{2} \mathbbm{E} \int_0^t \tau_{\ell} (V_{s -}) \mathd [\lambda
       Y, \lambda Y]^c_s + (\tau_{\ell} (V_{s -}))^{- 1} \mathd [\lambda^{- 1}
       Z, \lambda^{- 1} Z]^c_s \\
       && \hspace{2em}
       +\mathbbm{E} \sum_{0 < s \leqslant t} | \Delta Y_s |  | \Delta Z_s |\\
       & \geqslant& \mathbbm{E} \int_0^t  | \mathd [\lambda Y,
       \lambda^{- 1} Z]^c_s | +\mathbbm{E} \sum_{0 < s \leqslant t} | \Delta
       Y_s |  | \Delta Z_s |\\
       & \geqslant & \mathbbm{E} \int_0^t  | \mathd [Y, Z]^c_s |
       +\mathbbm{E} \sum_{0 < s \leqslant t} | \Delta Y_s |  | \Delta Z_s |\\
       & \geqslant & \mathbbm{E} \int_0^t  | \mathd [Y, Z]_s |
     \end{eqnarray*} 
  where all integrals and sums converge. Hence for all $\lambda > 0$
  \[ Q^{\mathcal{F}}_2 [w] (1 + \ell) (\lambda^2 \mathbbm{E}F_t + \lambda^{-
     2} \mathbbm{E}G_t + 2 a^2 \varepsilon^{- 1}) \gtrsim \mathbbm{E} \int_0^t
     | \mathd [Y, Z]_s | . \]
  We let now successively $\ell \rightarrow 0$ then $a \rightarrow 0$.
  Choosing specific $\lambda^2 = (\mathbbm{E}G_t)^{1 / 2} (\mathbbm{E}F_t)^{- 1 / 2}$,
  we can assume $\lambda > 0$ (otherwise the claim is trivial), we have
  \[ \mathbbm{E} \int_0^t | \mathd [Y, Z]_s | \lesssim Q^{\mathcal{F}}_2 [w] 
     (\mathbbm{E}F_t)^{1 / 2} (\mathbbm{E}G_t)^{1 / 2} \lesssim
     Q^{\mathcal{F}}_2 [w]  \| X \|_{2, w^{\varepsilon}}  \| Z \|_{2,
     u^{\varepsilon}} . \]
  The inequality above remains valid in the limit $t \rightarrow \infty$.
  Since the left--hand--side does not depend on the truncation of the weight,
  it remains to observe that
  \[ \lim_{\varepsilon \rightarrow 0} \| X \|_{2, w^{\varepsilon}} = \| X
     \|_{2, w} \hspace{1em} \tmop{and} \hspace{2em} \lim_{\varepsilon
     \rightarrow 0} \| Z \|_{2, u^{\varepsilon}} = \| Z \|_{2, u} . \]
  Indeed, since $X \in L^2 (\Omega ; \mathd \mathbbm{P}) \cap L^2 (\Omega ;
  \mathd \mathbbm{P}^w)$, we have for all $0 < \varepsilon < 1$, a.s.
  $X_{\infty}^2 w_{\infty}^{\varepsilon} \leqslant X_{\infty}^2 + X_{\infty}^2
  w_{\infty}$ and the limits above are a consequence of the dominated
  convergence theorem. The same reasoning applied to $Z$ completes the proof
  of the bilinear embedding.
\end{proof*}

\

\begin{proof*}{Proof of Theorem \ref{T: main result}}
  {\dueto{differential subordination under change of law}}The proof of the
  main result is now straightforward since the Proposition above allows us to
  estimate, for any test function $Z_{\infty} \in L^2 (\Omega, \mathd
  \mathbbm{P}) \cap L^2 (\Omega, \mathd \mathbbm{P}^u)$,
  \begin{eqnarray*}
    | (Y_{\infty}, Z_{\infty}) | & = & \left| \int_0^{\infty} \mathd [Y, Z]_s
    \right| \leqslant \int_0^{\infty} | \mathd [Y, Z]_s | \lesssim
    Q^{\mathcal{F}}_2 [w]  \| X \|_{2, w}  \| Z \|_{2, u},
  \end{eqnarray*}
  that is exactly
  \[ \| Y \|_{2, w} \lesssim Q^{\mathcal{F}}_2 [w]  \| X \|_{2, w} . \]
  This concludes the proof of Theorem \ref{T: main result}.
\end{proof*}

\

\section{Sharpness and applications}

\

\subsection{Sharpness}

\

Sharpness means that he linear power in the martingale $A_2$ characteristic cannot be improved. 

\subsubsection{discrete time}

That the result is sharp in the dyadic, discrete--in--time filtration case is well known and 
follows from the sharpness of the linear estimate for the dyadic square function in this setting (see
{\cite{HukTreVol2000a}} for an explicit calculation). Notice that the norm of
the square function is no larger than that of a predictable dyadic multiplier
- given the dyadic square function is obtained by taking expectation of a
$\sigma = \pm 1$ predictable multiplier $T_{\sigma}$. Indeed $S f^2 (t)
=\mathbbm{E} | T_{\sigma} f (t) |^2$, see for example {\cite{PetPot2002a}}.

\subsubsection{continuous time}

To see an example with continuous--in--time filtration, see
{\cite{DomPetWit2015a}} for details, we briefly summarize the flow of the
argument. Let $f (x)$ a compactly supported integrable function, \ defined on
$\mathbbm{R}$ and $\tilde{f} (z) = \tilde{f} (x, y)$ its harmonic extension to
the upper half space. Let $W_t = (x_t, y_t)$ be background noise (see
{\cite{GunVar1979}}), that is in a limiting sense a two-dimensional Brownian motion
starting at $\infty$ and arriving on the $x$ axis. Then the martingales
\[M_t^{\tilde{f}} = \tilde{f} (W_t) {\text{ and }} M_t^{\widetilde{H f}} = \widetilde{H
f} (W_t)\] ($H$ the Hilbert transform) are a pair of differentially subordinate
martingales that cannot have sublinear growth in weighted space with respect
to the $A_2$ characteristic of the induced filtration. To see this, use the
formula by {\tmname{Gundy--Varopoulos}} {\cite{GunVar1979}} restricted to the
Hilbert transform. By Cauchy--Riemann relations one passes to a martingale
representation that does not require conditioning by arrival, such as written
in the Riesz transform case in {\cite{GunVar1979}}. Then the authors in
{\cite{DomPetWit2015a}} borrowed the explicit examples that show the correct
growth of the Hilbert transform using the Poisson characteristic for $1 < p <
2$ and passed to $p = 2$ through an extrapolation argument using the
martingale setting through an argument by contradiction. Further, it is easy
to see that the deterministic Poisson $A_2$ characteristic and the martingale
$A_2$ characteristic driven by background noise are comparable.

\

\subsection{Applications}

\

\subsubsection{Discrete--in--time predictable multipliers}

The Bellman function in this paper and in particular its one-leg convexity can
give a direct proof of the results in {\cite{Lac2015a}} and
{\cite{ThiTreVol2015a}}, a weighted estimate for predictable multipliers in
the case of discrete in time filtrations. \

\subsubsection{Dimension--free weighted bounds on discrete operators}

Through the recent stochastic integral formula for second order Riesz
transforms {\cite{ArcDomPet2015b}} on compact multiply-connected Lie groups
$\mathbbm{G}$, our result gives dimension--free weighted $L^2$ estimates in
this setting too, using the semi-discrete heat characteristic of the weight.
The second order Riesz transforms take the form
\[ R^2_{\alpha} = \sum_i \alpha_i R_i^2 + \sum_{j, k} \alpha_{i j} R_{i j}^2,
\]
where the first diagonal sum are second order Riesz transforms in discrete
directions of the space and the second sum are continuous second order Riesz
transforms on the connected part, see {\cite{ArcDomPet2015b}} for more precise
definitions. The process considered is deterministic in one variable and is
Brownian in continuous directions together with a compound Poisson jump
process in the other, discontinuous directions. It was proved in
{\cite{ArcDomPet2015b}} that $R^2_{\alpha} f (z)$ can be written as the
conditional expectation $\mathbbm{E} (M_0^{\alpha, f} | Z_0 = z)$ where
$M^{\alpha, f}_t$ is a martingale transform of $M_t^f$ associated to $f$ and
$Z_t$ a suitable random walk. One obtains the estimate
\begin{equation}
  \| R_{\alpha} f \|_{L^2 (w)} \lesssim Q^{}_2 (w) \| f \|_{L^2 (w)}
  \label{estimate-weight}
\end{equation}
with implied constant independent of dimension and $Q_2 (w)$ the semi discrete
heat characteristic. An important special case are the second order discrete
Riesz transforms on products of integers. Notice that both the
continuous--in--time filtrations and the consideration of jump processes are
important to get this estimate.

It is also possible to get a deterministic proof of this application
(\ref{estimate-weight}), using the Bellman function we construct in this paper
in combination with part of the proof strategy in {\cite{DomPet2014a}}. Notice
though, that the trick used in {\cite{DomPet2014a}} to overcome the difficulty
of the jumps, does not work in the weighted setting, due to non-convexity of
the domain of the Bellman function. For a deterministic Bellman proof to give
the weighted estimate (\ref{estimate-weight}), it is instrumental to have the
one-leg convexity property we proved here.

\subsubsection{Probabilistic proof for estimate of the weighted Beurling
operator}

Our result gives a probabilistic proof of the weighted estimate for the
Beurling--Ahlfors transform that solved a famous borderline regularity problem
in {\cite{PetVol2002a}} previously proved by Bellman functions and other
means. To see this, one invokes the stochastic integral identity formula
{\cite{BanJan2008a}} for the Beurling-Ahlfors operator using heat flow
martingales. The comparability of heat flow $A_2$ characteristic and
martingale characteristic obtained when using the filtration in
{\cite{BanJan2008a}} is not hard to see. In turn, in {\cite{PetVol2002a}} it
was seen that the heat flow characteristic compares linearly to the classical
characteristic. The standard extrapolation result for sublinear operators in
{\cite{DraGraPerPet2005a}} gives the sharp weighted result in $L^p$.

\subsubsection{Dimension--free weighted bound, Riemannian setting}

{\tmname{Dahmani}} {\cite{Dah2016a}} used the continuous properties of the Bellman
function constructed in this paper to prove a dimensionless weighted bound for
the Bakry Riesz vector. Her result gives an optimal estimate in terms of the
Poisson characteristic. She considers a large class of manifolds with
non--negative Bakry-Emery curvature, such as for example the Gauss space. The
explicit expression of the Bellman function is essential to her argument. 

\

\

{\noindent}\begin{tabular}{l}
  \hline
  \quad
\end{tabular}

\


\begin{thebibliography}{10}
  \bibitem[1]{ArcDomPet2015b}Nicola Arcozzi, Komla Domelevo, and Stefanie
  Petermichl. {\newblock}Second order Riesz transforms on multiply-connected
  Lie groups and processes with jumps. {\newblock}\tmtextit{ArXiv e-prints},
  07 2015.
  
  \bibitem[2]{BanOse2014b}Rodrigo Ba{\~n}uelos and Adam Osekowski.
  {\newblock}On the Bellman function of Nazarov, Treil and Volberg.
  {\newblock}\tmtextit{Math. Z.}, 278(1-2):385--399, 2014.
  
  \bibitem[3]{BanOse2016a}Rodrigo Ba{\~n}uelos and Adam Osekowski.
  {\newblock}Sharp weighted $L^2$ inequalities for square functions.
  {\newblock}\tmtextit{ArXiv e-prints}, 03 2016.
  
  \bibitem[4]{BanJan2008a}Rodrigo Ba{\~n}uelos and Prabhu Janakiraman.
  {\newblock}$L^p$-bounds for the Beurling-Ahlfors transform.
  {\newblock}\tmtextit{Trans. Amer. Math. Soc.}, 360(7):3603--3612, 2008.
  
  \bibitem[5]{BerFrePet2015a}Fr{\'e}d{\'e}ric Bernicot, Dorothee Frey, and
  Stefanie Petermichl. {\newblock}Sharp weighted norm estimates beyond
  Calder{\'o}n-Zygmund theory. {\newblock}\tmtextit{ArXiv e-prints}, 10 2015.
  
  \bibitem[6]{BonLep1979a}A.~Bonami and D.~Lepingle. {\newblock}Fonction
  maximale et variation quadratique des martingales en pr{\'e}sence d'un
  poids. {\newblock}S{\'e}minaire de probabilit{\'e}s XIII, Univ. Strasbourg
  1977/78, Lect. Notes Math. 721, 294-306 (1979)., 1979.
  
  \bibitem[7]{Bur1984a}D.~L. Burkholder. {\newblock}Boundary value problems
  and sharp inequalities for martingale transforms. {\newblock}\tmtextit{Ann.
  Probab.}, 12(3):647--702, 1984.
  
  \bibitem[8]{Dah2016a}Kamilia Dahmani. {\newblock}Sharp dimension-free
  weighted bounds for the bakry-riesz vector. {\newblock}\tmtextit{In
  preparation}, 2016.
  
  \bibitem[9]{DelMey1982a}Claude Dellacherie and Paul-Andr{\'e} Meyer.
  {\newblock}\tmtextit{Probabilities and potential. B}, volume~72 of
  \tmtextit{North-Holland Mathematics Studies}. {\newblock}North-Holland
  Publishing Co., Amsterdam, 1982. {\newblock}Theory of martingales,
  Translated from the French by J. P. Wilson.
  
  \bibitem[10]{DomPet2014a}Komla Domelevo and Stefanie Petermichl.
  {\newblock}Bilinear embeddings on graphs. {\newblock}\tmtextit{preprint},
  2014.
  
  \bibitem[11]{DomPet2014c}Komla Domelevo and Stefanie Petermichl.
  {\newblock}Sharp $L^p$ estimates for discrete second order Riesz transforms.
  {\newblock}\tmtextit{Adv. Math.}, 262:932--952, 2014.
  
  \bibitem[12]{DomPetWit2015a}Komla Domelevo, Stefanie Petermichl, and Janine
  Wittwer. {\newblock}A linear dimensionless bound for the weighted riesz
  vector. {\newblock}\tmtextit{ArXiv e-prints}, 01 2015.
  
  \bibitem[13]{DraGraPerPet2005a}Oliver Dragi{\v c}evi{\'c}, Loukas Grafakos,
  Mar$\acute{\text{{\i}}}$a~Cristina Pereyra, and Stefanie Petermichl.
  {\newblock}Extrapolation and sharp norm estimates for classical operators on
  weighted Lebesgue spaces. {\newblock}\tmtextit{Publ. Mat.}, 49(1):73--91,
  2005.
  
  \bibitem[14]{DraTreVol2008a}Oliver Dragi{\v c}evi{\'c}, Sergei Treil, and
  Alexander Volberg. {\newblock}A theorem about three quadratic forms.
  {\newblock}\tmtextit{Int. Math. Res. Not. IMRN}, pages Art. ID rnn 072, 9,
  2008.
  
  \bibitem[15]{GunVar1979}Richard~F. Gundy and Nicolas~Th. Varopoulos.
  {\newblock}Les transformations de Riesz et les int{\'e}grales stochastiques.
  {\newblock}\tmtextit{C. R. Acad. Sci. Paris S{\'e}r. A-B}, 289(1):A13--A16,
  1979.
  
  \bibitem[16]{HukTreVol2000a}S.~Hukovic, S.~Treil, and A.~Volberg.
  {\newblock}The Bellman functions and sharp weighted inequalities for square
  functions. {\newblock}In \tmtextit{Complex analysis, operators, and related
  topics}, volume 113 of \tmtextit{Oper. Theory Adv. Appl.}, pages 97--113.
  Birkh{\"a}user, Basel, 2000.
  
  \bibitem[17]{HunMucWhe1973a}Richard Hunt, Benjamin Muckenhoupt, and Richard
  Wheeden. {\newblock}Weighted norm inequalities for the conjugate function
  and Hilbert transform. {\newblock}\tmtextit{Trans. Am. Math. Soc.},
  176:227--251, 1973.
  
  \bibitem[18]{Hyt2012a}Tuomas~P. Hyt{\"o}nen. {\newblock}The sharp weighted
  bound for general Calder{\'o}n-Zygmund operators. {\newblock}\tmtextit{Ann.
  Math. (2)}, 175(3):1473--1506, 2012.
  
  \bibitem[19]{IzuKaz1977a}M.~Izumisawa and N.~Kazamaki. {\newblock}Weighted
  norm inequalities for martingales. {\newblock}\tmtextit{T{\^o}hoku Math. J.
  (2)}, 29(1):115--124, 1977.
  
  \bibitem[20]{Lac2015a}Michael~T. Lacey. {\newblock}An elementary proof of
  the $A_2$ bound. {\newblock}01 2015.
  
  \bibitem[21]{LacPetReg2010a}Michael~T. Lacey, Stefanie Petermichl, and
  Maria~Carmen Reguera. {\newblock}Sharp $A_2$ inequality for Haar shift
  operators. {\newblock}\tmtextit{Math. Ann.}, 348(1):127--141, 2010.
  
  \bibitem[22]{Ler2015a}Andrei~K. Lerner. {\newblock}On pointwise estimates
  involving sparse operators. {\newblock}\tmtextit{ArXiv e-prints}, 12 2015.
  
  \bibitem[23]{NazTreVol1999a}F.~Nazarov, S.~Treil, and A.~Volberg.
  {\newblock}The Bellman functions and two-weight inequalities for Haar
  multipliers. {\newblock}\tmtextit{J. Amer. Math. Soc.}, 12(4):909--928,
  1999.
  
  \bibitem[24]{Ose2016a}Adam Osekowski. {\newblock}Sharp $L^p$-bounds for the
  martingale maximal function. {\newblock}\tmtextit{to appear in Tohoku
  Mathematical Journal}, 2016.
  
  \bibitem[25]{Pet2007a}S.~Petermichl. {\newblock}The sharp bound for the
  Hilbert transform on weighted Lebesgue spaces in terms of the classical
  $A_p$ characteristic. {\newblock}\tmtextit{Amer. J. Math.},
  129(5):1355--1375, 2007.
  
  \bibitem[26]{PetPot2002a}S.~Petermichl and S.~Pott. {\newblock}An estimate
  for weighted Hilbert transform via square functions.
  {\newblock}\tmtextit{Trans. Amer. Math. Soc.}, 354(4):1699--1703
  (electronic), 2002.
  
  \bibitem[27]{PetVol2002a}Stefanie Petermichl and Alexander Volberg.
  {\newblock}Heating of the Ahlfors-Beurling operator: weakly quasiregular
  maps on the plane are quasiregular. {\newblock}\tmtextit{Duke Math. J.},
  112(2):281--305, 2002.
  
  \bibitem[28]{Pro2005a}Philip~E. Protter. {\newblock}\tmtextit{Stochastic
  integration and differential equations}, volume~21 of \tmtextit{Stochastic
  Modelling and Applied Probability}. {\newblock}Springer-Verlag, Berlin,
  2005. {\newblock}Second edition. Version 2.1, Corrected third printing.
  
  \bibitem[29]{RezVasVol2010a}Alexander Reznikov, Vasiliy Vasyunin, and
  Alexander Volberg. {\newblock}An observation: cut-off of the weight $w$ does
  not increase the $a_{p_1, p_2}$-"norm'' of $w$. {\newblock}\tmtextit{ArXiv
  e-prints}, 08 2010.
  
  \bibitem[30]{ThiTreVol2015a}Christoph Thiele, Sergei Treil, and Alexander
  Volberg. {\newblock}Weighted martingale multipliers in the non-homogeneous
  setting and outer measure spaces. {\newblock}\tmtextit{Adv. Math.},
  285:1155--1188, 2015.
  
  \bibitem[31]{Tre2013a}Sergei Treil. {\newblock}Sharp $A_2$ estimates of Haar
  shifts via Bellman function. {\newblock}In \tmtextit{Recent trends in
  analysis. Proceedings of the conference in honor of Nikolai Nikolski on the
  occasion of his 70th birthday, Bordeaux, France, August 31 -- September 2,
  2011}, pages 187--208. Bucharest: The Theta Foundation, 2013.
  
  \bibitem[32]{VasVol2010c}Vasily Vasyunin and Alexander Volberg.
  {\newblock}Burkholder's function via Monge-Amp{\`e}re equation.
  {\newblock}\tmtextit{Ill. J. Math.}, 54(4):1393--1428, 2010.
  
  \bibitem[33]{Wan1995a}Gang Wang. {\newblock}Differential subordination and
  strong differential subordination for continuous-time martingales and
  related sharp inequalities. {\newblock}\tmtextit{Ann. Probab.},
  23(2):522--551, 1995.
  
  \bibitem[34]{Wit2000a}Janine Wittwer. {\newblock}A sharp estimate on the
  norm of the martingale transform. {\newblock}\tmtextit{Math. Res. Lett.},
  7(1):1--12, 2000.
\end{thebibliography}
\end{document}